\documentclass[11pt]{amsart}
\usepackage{amsbsy,amssymb,amscd,amsfonts,latexsym,amstext,delarray,
amsmath, graphicx, epsfig}  \headsep=15pt

\newtheorem{thm}{Theorem}[section]
\newtheorem{prop}[thm]{Proposition}
\newtheorem{cor}[thm]{Corollary}
\newtheorem{lem}[thm]{Lemma}

\newtheorem{defn}[thm]{Definition}

\numberwithin{equation}{section}

\def\R{{\mathbb R}}

\def\H{{\mathbb H}}

\def\O{{\mathcal O}}

\def\Q{{\mathbb Q}}
\def\C{{\mathbb C}}
\def\Z{{\mathbb Z}}

\def\P{{\mathbb P}}

\def\PSL{{\rm PSL}}
\def\PGL{{\rm PGL}}

\def\Aut{{\rm Aut}}
\def\Hom{{\rm Hom}}
\def\Tr{{\rm Tr}}
\def\sTr{{\rm sTr}}

\def\cP{{\mathcal P}}
\def\cK{{\mathcal K}}
\def\cS{{\mathcal S}}
\def\cM{{\mathcal M}}

\def\cH{{\mathcal H}}

\def\cA{{\mathcal A}}
\def\cV{{\mathcal V}}
\def\cN{{\mathcal N}}
\def\T{{\mathcal T}}
\def\cT{{\mathcal T}}
\def\cL{{\mathcal L}}

\def\cD{{\mathcal D}}
\def\cG{{\mathcal G}}
\def\cR{{\mathcal R}}
\def\cF{{\mathcal F}}
\def\cU{{\mathcal U}}
\def\cB{{\mathcal B}}

\newcommand{\ie}{{\it i.e.\/}\ }
\newcommand{\eg}{{\it e.g.\/}\ }
\newcommand{\cf}{{\it cf.\/}\ }

\title{Noncommutative geometry on trees and buildings}
\author[G.Cornelissen]{Gunther Cornelissen}
\author[M.Marcolli]{Matilde Marcolli}
\author[K.Reihani]{Kamran Reihani}
\author[A.Vdovina]{Alina Vdovina}
\address{G.~Cornelissen: University of Utrecht, The Netherlands}
\email{cornelissen@math.uu.nl}
\address{M.~Marcolli: Max--Planck Institut f\"ur Mathematik,
Bonn, Germany}
\email{marcolli\@@mpim-bonn.mpg.de}
\address{K.~Reihani: University of Oslo, Norway }
\email{kamranr\@@math.uio.no}
\address{A.~Vdovina: University of Newcastle, UK}
\email{Alina.Vdovina@newcastle.ac.uk}

\begin{document}

\maketitle

\section{Introduction}

The notion of a spectral triple, introduced by Connes (\cf
\cite{Connes}, \cite{Connes2}, \cite{CoMo}) provides a powerful 
generalization of
Riemannian geometry to noncommutative spaces. It originates from the
observation that, on a smooth compact spin manifold, the infinitesimal
line element $ds$ can be expressed in terms of the inverse of the
classical Dirac operator $\cD$, so that the Riemannian geometry is
entirely encoded by the data $(\cA,\cH,\cD)$ of the algebra of smooth
functions (a dense subalgebra of the $C^*$-algebra of continuous
functions), the Hilbert space of square integrable spinor sections,
and the Dirac operator. For a noncommutative space, the geometry is
then defined in terms of a similar triple of data $(\cA,\cH,\cD)$,
where $\cA$ is a $C^*$-algebra, $\cH$ is a Hilbert space on which
$\cA$ acts by bounded operators, and $\cD$ is an unbounded self
adjoint operator on $\cH$ with compact resolvent $(\cD-z)^{-1}$ for
$z\notin \R$, and such that the commutators $[D,a]$ are bounded
operators for all $a$ in a dense subalgebra of $\cA$.

\smallskip

Consani and Marcolli constructed in \cite{CM} a noncommutative space
describing the geometry of the 
special fibers at the archimedean places of an arithmetic surface, in
the form of a spectral triple for the action of a Kleinian Schottky group on its
limit set. The motivation for the construction was a result of Yuri Manin
\cite{Man}, computing the Arakelov Green function of a compact Riemann
surface in terms of a Schottky uniformization, and the proposed
interpretation of the ``dual graph'' of the fiber at arithmetic
infinity of an arithmetic surface in terms of the tangle of bounded
geodesics in a hyperbolic handlebody having the Riemann surface as its
conformal boundary at infinity (\cite{Man}, \cf also \cite{CM}).

\smallskip 

There is a well known analogy between archimedean places of an 
arithmetic surface and non-archimedean places with maximally degenerate
fiber. The fibers over such non-archimedean places also admit a
Schottky uniformization, by a p-adic Schottky group.  
This analogy was considered in  \cite{CM1}, where the construction of
\cite{CM} was generalized to the case of p-adic
Schottky groups and Mumford curves.
The use of Schottky uniformization in \cite{Man} implies that the
result of Manin on $\infty$-adic Arakelov geometry and hyperbolic
geometry appear to be confined to the 2-dimensional case. Some results
in higher dimension were obtained, in the case of linear cycles in
projective spaces, by Annette Werner \cite{Wer}, in terms of the geometry of the
Bruhat-Tits building for ${\rm PGL}(n)$. 

\smallskip

Motivated by this circle of ideas, we pursue two main directions in this
paper. One is a refinement of the construction of spectral triples for
Mumford curves. The main result is that we can improve the theta summable
construction of \cite{CM} to a {\em finitely summable} construction, upon passing
to the stabilization of the graph $C^*$-algebra of the dual graph of the
special fiber of a Mumford curve. The main advantage of a finitely summable 
spectral triple is that it makes it possible to extract invariants of the
geometry through zeta functions and through the Connes--Moscovici local 
index formula of \cite{CoMo}. The construction of the finitely summable
triple is based on a modification for the non-unital case
of a construction of Antonescu and Christensen \cite{AntChris} 
of spectral triples for AF algebras.

\smallskip

The other main direction we consider in the paper is that of
generalizations of the spectral triples from the case of Mumford
curves and trees to some classes of Euclidean and hyperbolic 
2-dimensional buildings.

\section{Theta summable spectral triples}\label{Sthetasum} 

We begin by recalling briefly the construction of the theta summable
spectral triple for the Kleinian Schottky case of \cite{CM} and give a
general formulation for group actions on trees, which also includes
the case of Mumford curves of \cite{CM1}. 

\medskip
\subsection{Kleinian Schottky groups}\label{KleinSect} 
\hfill\medskip

In the construction of \cite{CM}, one considers a Kleinian Schottky
group $\Gamma \subset \PSL(2,\C)$ acting by isometries on real
hyperbolic 3-space $\H^3$. This extends to an action on
$\P^1(\C)=\partial \H^3$ by fractional linear transformations. A
Kleinian Schottky group is a finitely generated discrete subgroup of
$\PSL(2,\C)$, which is 
isomorphic to a free group in $g$ generators and such that all
elements are hyperbolic. The
limit set $\Lambda_\Gamma\subset \P^1(\C)$ is the set of accumulation
points of orbits of the $\Gamma$-action. The quotient
$X=\Omega_\Gamma/\Gamma$, for $\Omega_\Gamma=\P^1(\C)\smallsetminus
\Lambda_\Gamma$, is a compact Riemann surface of genus $g$, which is
the conformal 
boundary at infinity of a hyperbolic handlebody of infinite volume
obtained as the quotient $\H^3/\Gamma$. The convex core
$H(\Lambda_\Gamma)/\Gamma$, where $H(\Lambda_\Gamma)\subset \H^3$ is
the convex hull of the limit set, is a region of finite volume and a
deformation retract of $\H^3/\Gamma$. For genus $g\geq 2$, the group
$\Gamma$ is non-elementary, namely the limit set $\Lambda_\Gamma$
consists of more than two points (it is in fact a totally disconnected
compact Hausdorff space -- a fractal in $\P^1(\C)$).

\smallskip

Consider the crossed product algebra $\cA=C(\Lambda_\Gamma)\rtimes
\Gamma$ of the action of the Schottky group on its limit set. This can
be identified with a Cuntz--Krieger algebra $O_A$ (\cf \cite{CuKrie}), 
where $A$ is the $2g\times 2g$ matrix with entries in $\{
0,1 \}$ associated to the  subshift of
finite type $\cS$ given by all admissible doubly infinite sequences in
the generators of the group $\Gamma$ and their inverses. 

\smallskip

The invertible shift $T$ on $\cS$ defines a noncommutative space
$C(\cS)\times_T \Z$. The corresponding homotopy quotient in the sense
of Baum--Connes is given by the quotient $\cS_T=(\cS\times \R)/\Z$. 
The cohomology $H^1(\cS_T,\C)$ is computed by an exact sequence
\begin{equation}\label{PV}
0 \to \C \to \cV \stackrel{\delta}{\to} \cV \to H^1(\cS_T,\C) \to 0,
\end{equation}
where $\cV$ is the infinite dimensional vector space
$\cV=C(\Lambda_\Gamma,\Z)\otimes \C$ of locally constant
functions on the limit set $\Lambda_\Gamma$. The coboundary $\delta$
is given by $\delta f = f- f\circ T$.

\smallskip

The space $\cV$ has a natural
filtration by finite dimensional vector spaces, where $\cV_n\subset
\cV$ is the space of locally constant
functions that only depend on the first $(n+1)$ coordinates. 
That is, for $\gamma\in \Gamma$ we let 
$\Lambda_\Gamma(\gamma)\subset \Lambda_\Gamma$ denote the set of
admissible infinite sequences $a_0, a_1, \ldots, a_n \ldots$ in the
generators and their inverses beginning with the word $\gamma$. Then
$\cV_n$ consists of continuous functions on $\Lambda_\Gamma$ that are
constant on each $\Lambda_\Gamma(\gamma)$ with $\gamma$ of length
$|\gamma|=n+1$ as an admissible word in the generators and their
inverses.  
The filtration of $\cV$ by the $\cV_n$ induces a filtration on the
cohomology 
$$ H^1(\cS_T,\C)=\varinjlim_n \cV_n/\delta \cV_{n-1}. $$
In \cite{CM} this cohomology was interpreted as a model for the
cohomology of the ``dual graph'' of the special fiber at 
infinity of an arithmetic surface.

\smallskip

The filtered vector space $\cV$ is a dense subspace of
$L^2(\Lambda_\Gamma, d\mu)$, where $d\mu$ is the Patterson--Sullivan
measure on the limit set $\Lambda_\Gamma$ (\cf \cite{Sull}), which
satisfies the scaling property 
\begin{equation}\label{PSmeas}
(\gamma^* d\mu)(x)= |\gamma^\prime (x)|^{\delta_H} \, d\mu(x), \ \ \
\forall \gamma\in \Gamma,
\end{equation}
with $\delta_H$ the Hausdorff dimension of $\Lambda_\Gamma$.

\smallskip

The grading operator $D=\sum_n\, n\, \hat\Pi_n$, where $\Pi_n$ is the
orthogonal projection onto $\cV_n$ and $\hat\Pi_n=\Pi_n-\Pi_{n-1}$,
is a densely defined self adjoint operator on $L^2(\Lambda_\Gamma,
d\mu)$ with compact resolvent. 

\smallskip

The Cuntz--Krieger algebra $O_A$ is the universal $C^*$-algebra
generated by $2g$ partial isometries $S_i$ subject to the relations 
\begin{equation} \label{CKrel} \sum_j S_j S_j^* =I \ \ \ \text{ and }
\ \ \ S_i^* S_i =\sum_j A_{ij} \, S_j
S_j^*, \end{equation}
where $A=(A_{ij})$ has entries $A_{ij}=1$ for
$|i-j|\neq g$, and $A_{ij}=0$ otherwise.
There is a faithful representation of the algebra $O_A$ on the algebra
of bounded operators $\cB(\cL)$ on the 
Hilbert space $\cL=L^2(\Lambda_\Gamma, d\mu)$, by setting
\begin{equation} (T_{\gamma^{-1}} f)(x)
:=|\gamma^\prime(x)|^{\delta_H/2}\, f(\gamma x),  \ \ \ \text{ and } \ \ \
(P_\gamma f)(x) := \chi_\gamma (x) f(x), \label{TandP}
\end{equation}
with $\chi_\gamma$ the characteristic function of the set
$\Lambda_\Gamma(\gamma)$, so that the $ S_i=\sum_j A_{ij} T_{\gamma_i}^*
P_{\gamma_j} $
are partial isometries satisfying \eqref{CKrel}.  

\smallskip

We give here a slightly modified version of the construction of \cite{CM}. 
Given an automorphism $U$ of the algebra $\cA$, and a representation
$\pi:\cA \to \cB(\cL)$ in the algebra of bounded operators of a
Hilbert space $\cL$, we consider the representation of $\cA$ on
$\cH=\cL\oplus \cL$ 
\begin{equation}\label{repU}
\pi_U(a)\,\,
(\xi,\zeta):=\left(\pi(a)\,\xi,\,\,\pi(U(a))\,\zeta\right). 
\end{equation}

\smallskip

Typically, for a Cuntz--Krieger algebra $\cA=O_A$, automorphisms of
the algebra can be
obtained from automorphisms of the corresponding space $\cS_A^+$ of
admissible right infinite sequences associated to the subshift of
finite type with matrix $A$. 
In fact, by \cite{Matsu}, there is a
homomorphism $\Aut(\cS_A^+) \to \Aut(O_A)$ that to 
an automorphism $u$ of $\cS_A^+$ assigns an automorphism $U$ of
$O_A$ which restricts to $u^*(f)=f\circ u^{-1}$ on the maximal
commutative subalgebra $C(\cS_A^+)$. By an automorphism of $\cS_A^+$
one denotes a homeomorphism $u$ of $\cS_A^+$ such that $T \circ u \circ
T^{-1}=u$, where $T$ is the one-sided shift on $\cS_A^+$.
In the case we are considering, we can
identify $\cS_A^+$ with the limit set $\Lambda_\Gamma$. 

\smallskip

A spectral triple $(\cA,\cH,\cD)$ for Kleinian Schottky groups is then
obtained as follows (\cf \cite{CM}).

\begin{prop}\label{S3KleinS}
Let $\cA=O_A$, let $\pi$ be the representation \eqref{TandP} on
$\cL=L^2(\Lambda_\Gamma,d\mu)$. Let $u$ be an automorphism of $\Lambda_\Gamma$ 
and let $\pi_U$ be the representation \eqref{repU} of $\cA$
on $\cH=\cL\oplus \cL$, for the induced automorphism $U$ of $\cA$.
Let $F$ be the linear involution that exchanges the two copies of $\cL$ and 
$\cD=FD$, with $D=\sum_n n \hat\Pi_n$. Then, for $\delta_H<1$, the
data $(\cA,\cH,\cD)$ define a spectral triple.
\end{prop}

\proof The result follows by showing that the commutators
$[\cD,a]$, for $a$ in a dense subalgebra of $O_A$, are bounded operators
on $\cH$. For that it is sufficient to estimate the norm of the
commutators $[D,S_i]$ and $[D,S_i^*]$ on $\cL$. An estimate is given
in \cite{CM} in terms of the Poincare' series of the Schottky group
(hence the $\delta_H<1$ condition). This takes care also of the
commutators with $\tilde S_i =U S_i$ and their adjoints. 
In fact, arguing as in \cite{Matsu} we see that the
$\tilde S_i=\sum_{j} S_{j} P_{\chi_{ij}}$, where $P_{\chi_{ij}}$
denotes the orthoginal projection associated to the 
characteristic function $\chi_{ij}$ of the set $E_{ij}$ of points
$\omega\in \Lambda_\Gamma$ such that the infinite word
$u \gamma_i u^{-1} \omega$ is admissible (\ie it defines a point in
$\Lambda_\Gamma$) and has $\gamma_j$ as first letter. In particular,
this implies that the automorphism $U$ preserves
the dense subalgebra of $\cA=O_A$ generated algebraically by the
partial isometries $S_i$ and their adjoints $S_i^*$, so that we still
have $[D,\tilde S_i]$ and $[D,\tilde S_i^*]$ bounded. 
\endproof

\smallskip

The spectral triple of Proposition \ref{S3KleinS} is $\theta$-summable
since the dimension of the eigenspaces of $D$ grows like
$2g(2g-2)^{n-1}(2g-2)$ so that we have
\begin{equation}\label{thetasum}
 \Tr(e^{-tD^2})<\infty \ \ \ \  \forall t>0. 
\end{equation}
The fact that the spectral triple is 
not finitely summable (which can be easily seen from the rate of
growth of the dimensions of the eigenspaces of $D$) falls under a
general result of Connes \cite{Connes2}, which shows that non-amenable
discrete groups (like Schottky groups) do not admit finitely summable
spectral triples. 

\smallskip

A reason for introducing the choice of an automorphism $U$ of the algebra
is in order to allow for a non-trivial $K$-homology class of the
spectral triple. It is often the case in specific cases of geometric 
interest that one has specific automorphisms available as part of the data.

\medskip
\subsection{Theta summable spectral triples for actions on trees}\label{TreesSect}
\hfill\medskip

We now present a similar construction of 
spectral triples in the case of actions on
trees, which refines the construction given in \cite{CM1} for the
case of Mumford curves. We show here that the construction
indeed follows very closely the case of Kleinian Schottky groups.

\smallskip

The results of Lubotzky \cite{Lu} show that the geometry of actions on
trees by finitely generated, torsion free, discrete subgroups of
automorphisms is in many ways analogous to that of Kleinian
Schottky groups. In the same philosophy, we shall see that the
construction of the spectral triple recalled in Section \ref{KleinSect}
also extends to this case. (This
in particular includes the case of Mumford curves.) The main results
that we need for this constructions are provided in \cite{Coo},
\cite{HeHu}, and \cite{Lu}.  

\smallskip

Let $\T$ be a locally finite tree, with $\T^0$ the set of vertices and
$\T^1$ the set of edges, and let $\Gamma \subset \Aut(\T)$ be a finitely
generated discrete subgroup. 
A path in $\T$ is a sequence $v_0, v_1, \ldots, v_n \ldots$ with
$v_i\in \T^0$, where $v_i$ and $v_{i+1}$ are adjacent and there are
no cancellations (namely $v_i\neq v_{i+2}$). The set of ends $\partial
\T$ is the set of equivalence classes of paths in $\T$, where
equivalent means having infinitely many $v_i$'s in common. A geodesic
in $\T$ is a doubly infinite path, namely a sequence $\ldots v_{-m},
\ldots, v_{-1}, v_0, v_1, \ldots, v_n \ldots$ with $v_i$ and
$v_{i+1}$ adjacent and $v_i\neq v_{i+2}$.
A distance function on the tree is obtained by assigning distance one 
to any pair of adjacent vertices.

\smallskip

The action of $\Gamma \subset \Aut(\T)$ extends to an action on
$\overline{\T}=\T \cup \partial \T$. The limit set
$\Lambda_\Gamma\subset \partial \T$ is
the set of accumulation points of $\Gamma$-orbits on
$\T$. 
Let $H(\Lambda_\Gamma)$ be the geodesic hull of the limit set, namely
the set of geodesics in $\T$ with both ends on $\Lambda_\Gamma$. It is
a closed $\Gamma$-invariant subset of $\T$, and
$H(\Lambda_\Gamma)/\Gamma$ is the convex core of $\T/\Gamma$. 

\smallskip

Lubotzky showed in \cite{Lu} that, if a finitely generated discrete
subgroup $\Gamma\subset \Aut(\T)$ is torsion free, then it is a
Schottky group, in the sense that $\Gamma$ is isomorphic to a free
group and every element is hyperbolic. 
In this case, the convex core $H(\Lambda_\Gamma)/\Gamma$ is a finite
graph. 

\smallskip

We first show how to extend naturally the construction of the
spectral triple $(\cA, \cH,\cD)$ of Proposition \ref{S3KleinS} to actions
on trees.

\smallskip

We will use essentially the fact that, by the result
of Coornaert \cite{Coo}, there is an analog for $\Gamma \subset
\Aut(\T)$ of the Patterson--Sullivan measure on the limit set of a
Kleinian Schottky group.
We follow the notation of \cite{HeHu} and assign to 
a hyperbolic element $\gamma\in \Aut(\T)$, the expression
\begin{equation}\label{horodist}
 (v,\gamma^{-1} v, x) = d(v,u) - d(\gamma^{-1} v,u), 
\end{equation}
where $v\in \T^0$ is a base point, $x$ is a point on the boundary
$\partial \T$, and $u$ is any vertex in the intersection of the two
paths from $v$ to $x$ and from $\gamma^{-1} v$ to $x$. The
distance $d(v,w)$ is the length of the geodesic arc in $\T$ connecting
$v,w\in \T^0$. The horospheric distance \eqref{horodist} does not
depend on $u$ and one defines
\begin{equation}\label{gammaprime}
\gamma^\prime_v(x)= e^{(v,\gamma^{-1} v, x)}, \ \ \ \text{ for } x\in
\partial \T.
\end{equation}
Then by \cite{Coo}, for $\delta_H$ the critical exponent
of the Poincar\'e series, there exists a normalized measure on
$\partial \T$ with support on $\Lambda_\Gamma$, satisfying
\begin{equation}\label{scalemeasure}
(\gamma^* d\mu_v)(x) = (\gamma^\prime_v(x))^{\delta_H} \, d\mu_v(x), \
\  \forall \gamma \in \Gamma.
\end{equation}
The Hausdorff dimension of $\Lambda_\Gamma$ is
equal to the critical exponent $\delta_H$.

\smallskip

As in the case of the Kleinian Schottky group, we can then consider
the Hilbert space $\cL=L^2(\Lambda_\Gamma, d\mu)$, with respect to the
$\Gamma$-conformal measure \eqref{scalemeasure}. The dense subspace of
locally constant functions $\cV=C(\Lambda_\Gamma,\Z)\otimes \C$ has a
filtration by $\cV_n$ defined as in the case of the Kleinian Schottky
group, by taking functions that depend only on the first
$(n+1)$-coordinates. Here we use an identification of
$\Lambda_\Gamma$ with admissible infinite sequences $a_0 a_1\ldots a_n
\ldots$ in the generators of $\Gamma$ and their inverses. Such
identification is determined by the choice of the base point $v\in \T^0$
and the identification $\Lambda_\Gamma = \overline{\Gamma v}\cap
\partial \T$.  

\smallskip

\begin{prop}\label{S3treeS} 
Let $\T$ be a locally finite tree and $\Gamma \subset \Aut(\T)$ be a
torsion free finitely generated discrete subgroup with $\Lambda_\Gamma
\subset \partial \T$ its limit set. Let $u$ an automorphism of
$\Lambda_\Gamma$ and $U$ the induced automorphism of the $C^*$-algebra 
$\cA=C(\Lambda_\Gamma)\rtimes \Gamma$.  
Then the data $(\cA, \cH,\cD)$ as in Proposition
\ref{S3KleinS} define a $\theta$-summable spectral triple.
\end{prop}

\proof For $\Gamma\subset \Aut(\T)$ a torsion free
finitely generated discrete subgroup (hence a Schottky group by
\cite{Lu}), let $\{ \gamma_i \}_{i=1}^{g}$ be a set of generators 
and let $\gamma_{i+g}=\gamma_i^{-1}$. 

\smallskip

The representation of the algebra $\cA=C(\Lambda_\Gamma)\rtimes
\Gamma$ on the Hilbert space $\cL=L^2(\Lambda_\Gamma, d\mu)$ is then
given as in \eqref{repU}, with $U$ induced by an automorphism $u$ of
$\Lambda_\Gamma$ and the representation $\pi$ on $\cB(\cL)$ defined as
in \eqref{TandP} by setting 
\begin{equation}\label{TandP2}
(T_{\gamma^{-1}} f)(x)
:=(\gamma^\prime_v(x))^{\delta_H/2}\, f(\gamma x),  \ \ \ \text{ and }
\ \ \ (P_\gamma f)(x) := \chi_\gamma (x) f(x),
\end{equation}
where $\chi_\gamma$ is the characteristic function of the subset
$\Lambda_\Gamma(\gamma) \subset \Lambda_\Gamma$.

\smallskip

As in the proof given in \cite{CM} of the result for Kleinian Schottky
groups of Proposition \ref{S3KleinS}, it is then sufficient to prove that,
for $i=1,\ldots, 2g$, the commutators $[D,S_i]$ and $[D,S_i^*]$ are
bounded. Here $S_i=\sum_j A_{ij} T_{\gamma_i}^*
P_{\gamma_j}$ is the operator
$(S_i\, f)(x)=(1-\chi_{\gamma_i^{-1}}(x))
(\gamma_i^\prime(x))^{\delta_H/2} \, f(\gamma_i x)$,
while the adjoint $S^*_i$ acts as
$(S_i^* f)(x)=\chi_{\gamma_i}(x) (\gamma_i^\prime
(Tx))^{-\delta_H/2} f(Tx)$,
where $T$ is the one sided shift on $\Lambda_\Gamma$ satisfying
$T|_{\Lambda_\Gamma(\gamma_i)}(x)=\gamma_i^{-1} (x)$.

\smallskip

We use the fact (\cite{Coo}, \cite{HeHu}) that the function
$\gamma^\prime_v(x)$ is locally constant on
$\partial \T$, for any chosen $v\in \T^0$. This means that 
$$ \Pi_{k_i} \, ( \gamma_i^\prime) = 
\gamma_i^\prime , $$
for some $k_i>0$. This implies that we have
$$ S_i : \cV_{n+1} \to \cV_n, \ \ \ \  S_i^* : \cV_n \to
\cV_{n+1}, \ \ \ \ \forall n\geq k_i. $$
Thus the commutator $[D,S_i]$ can be written as 
\begin{equation}\label{commDSi}
 [D,S_i] =-S_i(1-\Pi_0) + \sum_{k=0}^{k_i-1} (S_i \Pi_{k+1} - \Pi_k
S_i), 
\end{equation}
since $S_i \Pi_{k+1}=\Pi_k S_i$ for $k\geq k_i$. Thus $[D,S_i]$ is a
bounded operator. The argument for $[D,S^*_i]$ is analogous. 
Thus the commutators $[D,a]$ are bounded for all elements $a$
in the dense involutive subalgebra of $C(\Lambda_\Gamma)\rtimes
\Gamma$ generated algebraically by the $S_i$. 

\smallskip

We use again, as in Proposition \ref{S3KleinS} above, the fact that
the automorphism $U$ of $\cA$ induced by the automorphim $u$ of 
$\Lambda_\Gamma$ preserves the subalgebra generated algebraically by
the $S_i$ and $S_i^*$ so that we have bounded commutators with elements 
of this dense subalgebra in the representation twisted by $U$.
The rest of the argument is then analogous to \cite{CM} and shows that we
obtain a spectral triple. 

\endproof

\smallskip

This construction can be refined by working with the graph
$H(\Lambda_\Gamma)/\Gamma$, instead of directly with the limit set
$\Lambda_\Gamma$. 
The hull $H(\Lambda_\Gamma)$ consists of all the axes $L(\gamma)$ 
of $\gamma\in\Gamma$. These are the geodesics in $\T$ connecting the
fixed points $z^-(\gamma)$ and $z^+(\gamma)$ in $\Lambda_\Gamma$. Let
$A$ be the {\em directed edge matrix} of the finite graph
$\cG=H(\Lambda_\Gamma)/\Gamma$. This is a $2 \# \cG^1\times 2 \#
\cG^1$ matrix with entries in $\{0,1\}$, such that $A_{ee'}=1$ is
$ee'$ is an admissible path in $\cG$ where $e,e'$ are edges with
either possible orientation, and $A_{ee'}=0$ otherwise. The set
$\cS_A^+$ of infinite admissible words, with the admissibility
condition specified by the matrix $A$, describes the set of 
infinite walks on the tree $H(\Lambda_\Gamma)$ starting at any vertex
of a given fundamental domain for the action of $\Gamma$. This gives an
identification of $\cS_A^+$ with a union of $n$ copies of
$\Lambda_\Gamma$, for $n=\# \cG^0$, each corresponding to walks
starting at a given vertex in $\cG$. Using this identification, we can
define the Hilbert space $\cL_A=L^2(\cS_A^+,d\mu)$, with the measure induced
by the Patterson--Sullivan measure \eqref{scalemeasure} on
$\Lambda_\Gamma$. This identification also induces on $\cL_A$ a
filtration, for which we can consider the associated grading $D$.
As in the previous cases, we have a representation
of the algebra $\cA=C(\cS_A^+)\rtimes \Gamma$ on the Hilbert space
$\cL_A$. The algebra $\cA$ is Morita equivalent to a graph
$C^*$-algebra of $H(\Lambda_\Gamma)/\Gamma$. As before, we can
consider an automorphism of the Cantor set $\cS_A^+$ and a
corresponding automorphism $U$ of the algebra $\cA$ and a
representation on $\cH_A=\cL_A\oplus\cL_A$, where the $[F,\pi_U(a)]$ 
are compacts (\cf Lemma 4.6, \cite{Matsu}).
We can then extend the result of
Proposition \ref{S3treeS} to $\theta$-summable spectral triples of the form
$(C(\cS_A^+)\rtimes \Gamma, \cH_A,\cD=FD)$.

\smallskip

The original spectral triple of Proposition \ref{S3treeS} corresponds 
to considering the graph $\T_\Gamma/\Gamma$ instead of
$H(\Lambda_\Gamma)/\Gamma$, where $\T_\Gamma$ is the Cayley graph of
$\Gamma$ and $H(\Lambda_\Gamma)$ is the smallest subtree of $\T$
containing all the axes of the elements of $\Gamma$, \ie the infinite
geodesics in $\T$ with endpoints $z_\gamma^\pm\in \Lambda_\Gamma$ the
fixed points of $\gamma$. In the case of Mumford curves, where $\T$ is
the Bruhat--Tits tree of $\PGL(2,K)$, with $K$ a finite extension of
$\Q_p$, the graph $\T_\Gamma/\Gamma$ gives the dual graph of the
specialization over the ring of integers $\O\subset K$ of the
algebraic curve $C$ holomorphically isomorphic to
$X=\Omega_\Gamma/\Gamma$, while $H(\Lambda)/\Gamma$ is the dual
graph of the minimal smooth model of $C$ over $\O$, \cf \cite{Mum}.  

\smallskip

The difference between the graphs
$\T_\Gamma/\Gamma$ and $H(\Lambda)/\Gamma$ may be used to define a
``code building'' procedure of the type considered in \S 5 of
\cite{LiMa} to produce codes from Mumford curves. We shall not
deal with this aspect in the present paper.

\medskip

The above construction will also work if the group $\Gamma$ is replaced by a general finitely generated discrete subgroup $N$ of $\PGL(2,K)$; on the 
side of 
Mumford curves this corresponds to orbifold uniformization, or Mumford curves with automorphisms, cf.\ \cite{CoKa}. Such $N$ always has a finite index normal free subgroup. Conversely, given such a free group $\Gamma$,
a finitely generated discrete subgroup $N\subset
\PGL(2,K)$ contained in the normalizer $N(\Gamma)$ of $\Gamma$ in
$\PGL(2,K)$ determines a finite group $G=N/\Gamma \hookrightarrow
\Aut(X)$ of automorphisms of the Mumford curve
$X=\Omega_\Gamma/\Gamma$, since $\Aut(X)=N(\Gamma)/\Gamma$. 
Thus, it becomes relevant to study the
equivariant deformation problem, of how these data can be deformed to
another curve of the same genus with an action of the same group (\cf
\cite{CoKa}). For $\Gamma \subset N$ as above, the finite group
$G=N/\Gamma$ acts on $\cG= H(\Lambda_\Gamma)/\Gamma$ with quotient 
the finite graph $\cG_N=H(\Lambda_\Gamma)/N$.
Let $\rho_0: N\hookrightarrow \Aut(\T)$ denote the inclusion of
$N\subset N(\Gamma)\subset  \Aut(\T)$. Let $\Hom^*(N,\Aut(\T))$ denote
the subset of $\Hom(N,\Aut(\T))$ of injective homomorphisms with
discrete image. This governs the equivariant deformations. 
There is an open neighborhood of $\rho_0$ in
the space $\Hom^*(N,\Aut(\T))$ of the form
\begin{equation}\label{neighbHom}
\cU(\rho_0)=\{ \rho \in \Hom^*(N,\Aut(\T))| \,\,
\rho(\gamma_i)(x_i)=y_i, \,\, \rho(\gamma_i)(y_i)=\gamma_i(y_i) \}.
\end{equation}
Here $\{ \gamma_i \}_{i=1}^g$ is a set of generators for the Schottky
group $\Gamma \subset N$ and $x_i, y_i$ are points on the axes
$L(\gamma_i)$ that specify the Schottky data. As shown in \cite{HeHu}, 
the neighborhood $\cU(\rho_0)$ has the property that, for all
$\rho\in\cU(\rho_0)$, the group $\rho(\Gamma)$ is a Schottky group, of
finite index in $\rho(N)$, with the same Schottky data as the original
$\rho_0(\Gamma)$. In particular, in our setting, this shows that, 
by themselves, the spectral triples introduced in the previous 
section will not distinguish Mumford curves in the family  
$\rho(\Gamma)\subset \rho(N)$, with $\rho\in \cU(\rho_0)$.
However, one can implement in the construction the induced action
of the group $G$ on the $C^*$-algebra of the graph 
$\cG= H(\Lambda_\Gamma)/\Gamma$, or consider $G$-equivariant
spectral geometries for the Morita equivalent algebra
$C(\cS_A^+ \times G)\rtimes N$, with $N$ acting on the left on 
$G=N/\Gamma$. 

\medskip

In \cite{CM1} the local L-factor $L_v(H^1(X),s)=\det (1-Fr_v^*
N(v)^{-s}|H^1(\bar X,\Q_\ell)^{I_v})^{-1}$ of a Mumford curve was
recovered from the data $(\cA, \cH_A,\cD)$. Another possible direction
in which the construction of such spectral triples may be or
arithmetic significance is by associating to a Mumford curve some
cocycles in the entire cyclic cohomology of a smooth subalgebra of
$\cA$. 

\smallskip

Cyclic cohomology was introduced by Connes as a natural receptacle 
for the characters 
of finitely summable Fredholm modules. Similarly, in the theta
summable case that corresponds to ``infinite dimensional geometries'',
one can also define characters through the JLO cocycle 
$\varphi=(\varphi_{2n})$ (\cf \cite{Co94}) of the form
\begin{equation}\label{JLOeven}
\varphi_{2n}(a^0,\ldots, a^n)= \displaystyle{\int_{\substack{ s_i\geq
0 \\ \sum s_i=1}}} 
\sTr\left(a^0 e^{-s_0\cD^2} [\cD,a^1]e^{-s_1\cD^2}\cdots
[\cD,a^{2n}]e^{-s_{2n}\cD^2}\right) \,\, ds_0 \cdots ds_{2n} ,
\end{equation}
where $\sTr$ denotes the supertrace. These live naturally in the entire
cyclic cohomology of \cite{Co94}.

\section{Finitely summable spectral triples}\label{finsum3}

It is clear that the $\theta$-summable condition imposes a
strong limitation on how one may be able to apply tools from
noncommutative geometry to the
arithmetic context, most notably the local index formula of
\cite{CoMo}, which requires the finitely summable setting. 
There were already strong indications from the
original construction (\cf \cite{CM2}) that it should be possible to
obtain a finitely summable spectral triple associated to Mumford
curves. 

In fact the Cuntz--Krieger algebra $O_A$ has, up to stabilization
(\ie tensoring with compact operators) a second description as
a crossed product (\cf \cite{CuKrie}). Namely, one has 
an identification
\begin{equation} \label{AF-T} \overline{O_A} \cong
\overline{\mathcal F}_A \rtimes_T \Z, \end{equation}
where $\cF_A$ is an approximately finite dimensional (AF) algebra, \ie a
direct limit of finite dimensional algebras.
Here for a unital $C^*$-algebra $\cA$ we use the notation
$\overline{\cA} =\cA \otimes \cK$ with $\cK$ the algebra of
compact operators. The algebra $\overline{\cA}$ is no longer unital.

\smallskip

In the case of the algebra $O_A\cong C(\Lambda_\Gamma)\rtimes \Gamma$,
the hyperbolic growth of the 
Schottky group $\Gamma$ prevents one from constructing a finitely
summable Dirac operator. This is no longer the case for the
algebra $\cA_\Gamma:=\overline{\cF_A}\rtimes \Z$.
The fact that this algebra can be written as a crossed product by the
integers implies that, by Connes' result on hyperfiniteness
\cite{Connes2}, it may carry a finitely summable spectral triple. 

\smallskip

However, the fact of working with non-unital algebras forces one to
relax the axioms of finitely summable spectral triple to a suitable
``local'' version, as discussed in \cite{Moyal} in the important
example of Moyal planes. 

\smallskip

\begin{defn}\label{noncompS3}
Let $\cA$ be a non-unital $C^*$-algebra. A spectral triple
$(\cA,\cH,\cD)$ consists of the data of a representation
$\pi: \cA \to \cB(\cH)$ of the algebra as bounded operators on a
separable Hilbert space $\cH$, together with an unbounded self-adjoint
operator $\cD$ on $\cH$ such that the subalgebra
\begin{equation}\label{Ainfty}
\cA_\infty:=\{ a \in \cA\, \mid \, a{\rm Dom}\cD \subseteq {\rm Dom}\cD,\,\,
[\cD,a]\in \cB(\cH), \,\, a(1+\cD^2)^{-1}\in\cK(\cH)\}
\end{equation}
is dense in $\cA$.
\end{defn}

In \eqref{Ainfty} $\cK(\cH)$ denotes the ideal of compact operators.
In particular, if $\cD$ has compact resolvent,
then the last property of $\cA_\infty$ is automatically fulfilled.

\smallskip

In the case we are interested in, the AF algebra $\cF_A$ can be described
in terms of a groupoid 
$C^*$-algebra associated to the ``unstable manifold'' in the
Smale space $({\mathcal S},T)$. In fact, consider the algebra
${\mathcal O}_A^{alg}$ generated algebraically by the $S_i$ and
$S_i^*$ subject to the Cuntz--Krieger relations \eqref{CKrel}.
Elements in ${\mathcal O}_A^{alg}$ are linear
combinations of monomials $S_\mu S_\nu^*$, for multi-indices $\mu$,
$\nu$, \cf \cite{CuKrie}. The AF algebra ${\mathcal F}_A$ is generated
by elements $S_\mu S_\nu^*$ with $|\mu|=|\nu|$, and is filtered by
finite dimensional algebras ${\mathcal F}_{A,n}$ generated by elements
of the form $S_\mu P_i S_\nu^*$ with $|\mu|=|\nu|=n$ and
$P_i=S_iS_i^*$ the range projections, and embeddings determined by the
matrix $A$. 
The commutative algebra ${\rm C}(\Lambda_\Gamma)$ sits as
a subalgebra of ${\mathcal F}_A$ generated by all range projections
$S_\mu S_\mu^*$. The embedding is compatible with the filtration and
with the action of the shift $T$, which is implemented on ${\mathcal
F}_A$ by the transformation $a\mapsto \sum_i S_i\, a\, S_i^*$. (\cf
\cite{CuKrie}.)
The stabilization $\overline{\cF_A}$ is a non-unital AF algebra.

\smallskip

We use the description \eqref{AF-T} 
as the starting point for the construction of finitely summable
spectral triples. 

\smallskip

The following result provides a modification of Theorem 2.1
of \cite{AntChris}, where the AF-algebra is now not necessarily
unital, the sequence of eigenvalues need not 
be positive, and the growth condition required for the sequence 
of eigenvalues is finer. This will suffice for our purpose of
defining a finitely summable spectral triple.

\begin{thm}\label{AFnon1}
Let $\cA$ be a (not necessarily unital) AF-algebra, and let $p$ 
be a positive real number. 
Then there exists an unbounded Fredholm module $(\cH,\cD)$ over 
$\cA$ such that $(1+\cD^2)^{-p/2}\in \cL^1(\cH)$. In particular, 
$(\cA,\cH,\cD)$ is a $p$-summable odd spectral triple.
\end{thm}

\proof
We can write $\cA=\overline{{\cup}_{n=1}^\infty \cA_n}$, where each
$\cA_n$ is a finite dimensional $C^*$-algebra. 

Let $\tau$ be a state on $\cA$, that is, a continuous linear 
functional  $\tau: \cA\to \C$ of norm one, such that $\tau(a^*a)\geq 0$ 
for all $a\in \cA$.
The norm of a positive linear functional is the limit of its evaluation
at any approximate identity, so that, in the unital case a state is just
a positive linear functional with $\tau(1)=1$. We denote by 
$\cH=L^2(\cA,\tau)$ the
Hilbert space of the GNS representation defined by the state
$\tau$. Namely, $\cH$ is the Hilbert space completion of the 
quotient $\cA/\cN_\tau$, for $\cN_\tau$ is the closed left ideal  
$\cN_\tau=\{x \in \cA: \tau(x^*x)=0\}$, with respect to the inner product 
induced by $\tau(b^*a)$. 
In this non-unital case the cyclic vector $\xi$ in the GNS representation 
is obtained as the limit in the norm of $\cH$ of the classes of 
a given approximate identity for $\cA$.

Let $\eta$ denote the quotient map $\eta: \cA \to \cH$. 
We set $\cH_n=\eta(\cA_n)$.
These are finite dimensional subspaces of $\cH$ with $\dim \cH_n\leq \dim
\cA_n$ and $\cH_n\subset \cH_{n+1}$. We assume that the $\cH_n$ give a
filtration of $\cH$. 
Let $\Pi_n:\cH \rightarrow \cH_n$ be the orthogonal projection onto $\cH_n$ 
and put $\hat\Pi_n=\Pi_n-\Pi_{n-1}$, for $n \geq 2$, and $\hat\Pi_1
=\Pi_1$. 

Following the construction of \cite{AntChris},
we now show that we can choose a sequence $(\lambda_n)$ of {\em real} 
numbers such that the unbounded operator $\cD=\sum_{n=1}^\infty\lambda_n \hat\Pi_n$ 
defined on the dense subspace $\cup_{n=1}^\infty \cH_n$ of $\cH$ 
satisfies the $p$-summability condition. 

In fact, for any $n$ we may assume 
that $\cH_n \subsetneq \cH_{n+1}$, since otherwise we would have
$\hat\Pi_{n+1}=0$, and $\cA_n \subsetneq \cA_{n+1}$, so that 
we have $\dim \cA_{n+1}> \dim \cA_n$. This gives  
$\dim \cA_n \ge n$, since $\dim \cA_1 \ge 1$. Now, if we choose 
the eigenvalues $\lambda_n$ so that $|\lambda_n|\ge(\dim \cA_n)^q$,
for some $q>2/p$, we obtain an estimate of the form
\begin{align*}
\Tr(1+\cD^2)^{-p/2}&=\sum_{n=1}^\infty(1+|\lambda_n|^2)^{-p/2}
\dim E_{\lambda_n}\\&=\sum_{n=1}^\infty(1+|\lambda_n|^2)^{-p/2}
(\dim \cH_n - \dim \cH_{n-1})\\&\le \sum_{n=1}^\infty|\lambda_n|^{-p}
(\dim \cH_n) \le \sum_{n=1}^\infty|\lambda_n|^{-p}
(\dim \cA_n)  \\&\le\sum_{n=1}^\infty(\dim \cA_n)^{-p q}(\dim \cA_n) 
\le\sum_{n=1}^\infty \frac{1}{n^{pq-1}}< \infty~,
\end{align*}
where we set $\cA_0=\{0 \}$. Also notice that, for $m\geq n$ we have
$\cA_m(\cH_n)\subset \cH_n$, so that for $n > m$ and 
$a \in \cA_m$ we have $[\hat\Pi_n,a]=0$. Thus, we obtain 
$[\cD,a]=\sum_{n=1}^m \lambda_n [\hat\Pi_n,a]$. This shows that 
$[\cD,a]$ has a bounded closure on $\cH$. Moreover, if $a \in \cA_m$ 
then $a{\rm Dom}\cD \subseteq {\rm Dom}\cD$. Thus, the subalgebra 
$$\cA_\infty:=\{a \in \cA \mid a{\rm Dom}\cD \subseteq {\rm Dom}\cD~{\rm and}~
[\cD,a]~\textnormal{admits a bounded closure}\}$$ contains 
$\cup_{n=1}^\infty \cA_n$, hence it is dense in $\cA$.
\endproof

In fact, one can strenghten the result of the theorem as in Corollary
\ref{faithstate} below. However, we are interested in allowing for 
possibly non-faithful states, since we are interested in a state that
encodes the data of the uniformization through the Patterson--Sullivan 
measure on the limit set as in the theta summable construction.

\begin{cor}\label{faithstate}
In the construction of Theorem \ref{AFnon1} 
one can always choose the state $\tau$
to be faithful.
\end{cor}

\proof
It is enough to show that every separable $C^*$-algebra admits a faithful state.
Let $\cA$ be a separable $C^*$-algebra and let $\cA_+$ be its positive
cone. Choose a sequence $(a_m)$ which is dense in  $\cA_+$, and choose
a bounded approximate identity $(u_n)$ for $\cA$ in $\cB_1(\cA_+)$. Set 
$w_{n,m}:=u_n a_m u_n$ in $\cA_+$, and choose a state $\tau_{n,m}$ on 
$\cA$ such that $\tau(w_{n,m})=\|w_{n,m}\|$. Now define 
$$\tau:=\sum_{n,m \ge 1}^\infty \frac{\tau_{n,m}}{2^{n+m}}.$$
One can see that $\tau$ is a faithful state on $\cA$. 
\endproof

In the case of a faithful state, the Hilbert space $\cH$ is the
closure of $\cA$ in the inner product $\langle a,b\rangle =\tau(b^*a)$
and the filtration of $\cH$ is given by $\cH_n=\eta(\cA_n)=\cA_n$.

\smallskip

We then obtain a finitely summable spectral triple for the
algebra \eqref{AF-T} through the following construction.

\begin{thm}\label{S3cross}
Let $(\cA,H,D)$ be an odd spectral triple for the (not necessarily unital) 
$C^*$-algebra $\cA$, and assume that $\cD$ has compact
resolvent. Consider the crossed product  
$\cA \rtimes_\alpha \Z$, and assume that the dense subalgebra 
$\cA_\infty$ of the spectral triple contains a dense $\alpha$-invariant
subalgebra $\cA_\infty '$ such that, for any fixed $a \in \cA_\infty
'$ the sequence of bounded operator
$\{[\cD,\alpha^n(a)]\}_{n \in\Z}$ is uniformly bounded in the operator
norm. Let $(C(S^1),\ell^2(\Z),\partial)$ be the standard spectral triple
on the circle. 
Consider then the data 
\begin{equation}\label{crossS3}
(\cA\rtimes_\alpha\Z,\cH,\cD)
\end{equation}
where $\cH=\ell^2(\Z,H)\oplus\ell^2(\Z,H)$ is the direct sum of
two copies of the Hilbert space for the regular representation of 
$\cA\rtimes_\alpha\Z$ (as a reduced crossed product) and 
$$ \cD=\left(\begin{array}{cc}
0 & D^*_0\\
D_0 & 0\\
\end{array}\right),$$ 
for $D_0=D\otimes 1+i\otimes\partial$.
The data \eqref{crossS3} define an even spectral triple for
the algebra $\cA\rtimes_\alpha\Z$, with respect to the grading
on $\cH$ given by 
$$\gamma=\left(\begin{array}{cc}
1 & 0\\
0 & -1\\
\end{array}\right).$$
If $(\cA,H,D)$ is $(p,\infty)$-summable,
then $(\cA\rtimes_\alpha\Z,\cH,\cD)$ is 
$(p+1,\infty)$-summable.
\end{thm}
 
\proof Notice $\cD$ is in fact
the tensor product of $D$ and $\partial$ in the Baaj-Julg's picture
\cite{sB83} of Kasparov's external product in $K$-homology
$$ K^1(\cA)\times K^1(C(S^1))\rightarrow K^1(\cA \otimes C(S^1)). $$
In particular, the facts that $\cD$ is a selfadjoint operator with
compact resolvent and is finitely summable is already well known.

Let $V$ denote the regular representation of $C(S^1)=C^*(\Z)$
in $\ell^2(\Z)$, as part of the data $(C(S^1),\ell^2(\Z),
\partial)$. Put $\cL=\ell^2(\Z,H)$, so that $\cH
=\cL\oplus\cL$ and $\cA\rtimes_\alpha\Z$ acts 
diagonally on $\cH$. Then, for all 
$\xi\in\cL$, one has
\begin{equation}\label{cD0}
(\cD_0 \xi)(k)=D(\xi(k))-ik\xi(k), \ \ \ 
\forall k \in\Z.
\end{equation}
An element $a \otimes V^n \in C_c(\Z,\cA)\subset
\cA\rtimes_\alpha\Z$, for $a \in \cA'_\infty$ and $n \in\Z$,
is represented in $\cL$ by
$$((a \otimes V^n)\xi)(k)=\alpha^k(a)(\xi(n+k)), \ \ \ \forall k
\in\Z.$$ 
Using \eqref{cD0} one obtains 
$$[\mathcal{D},a \otimes V^n]=\left(\begin{array}{cc}
0 & T_1+inT_2\\
T_1-inT_2 & 0\\
\end{array}\right),$$ 
where
$T_i:\cL\rightarrow\cL$ are (unbounded) operators for
$i=1,2$ given by
$(T_1 \xi)(k)=[D,\alpha^k(a)](\xi(n+k))$
and
$(T_2 \xi)(k)=\alpha^k(a)(\xi(n+k))$, 
for $k \in\Z$. These satisfy the estimates
$\|T_1\|\le\sup_{k}\|[D,\alpha^k(a)]\|<\infty$, and 
$\|T_2\|\le\|a\|<\infty$. Therefore we see that 
$[\cD,a \otimes V^n]$ 
admits a bounded closure in $\cB(\cH)$.
\endproof

\smallskip

{}From an index theoretic perspective, it may be preferable to 
work with a finitely summable even spectral triples on the AF-algebra,
since it is the $K_0$-group of an AF algebra that carries all the
interesting information, while the $K_1$-group is trivial. It is
easy to modify the previous construction to accommodate this
case. As before, let $\cA=\overline{{\cup}_{n=1}^\infty \cA_n}$ be an
AF-algebra and let $\tau$ be a state on $\cA$. Consider the
Hilbert space $H=L^2(\cA,\tau)\oplus L^2(\cA,\tau)$. One can take
the diagonal action of $\cA$ on $H$ in the GNS representation,
although it is better to proceed as in the $\theta$-summable case and 
introduce a twisting on one of the copies of $L^2(\cA,\tau)$ by a
nontrivial automorphism of the algebra. 
Now one can choose  
a sequence of {\em complex} (not necessarily real) numbers $\lambda_n$ satisfying 
a suitable growth condition (\eg $|\lambda_n|\ge(\dim \cA_n)^q$, where $q>2/p$) as
the eigenvalues of an operator $D_0=\sum_{n=1}^\infty\lambda_n \hat\Pi_n$
on $L^2(\cA,\tau)$. One then consider the operator on $\cH$ defined by
$$D=\left(\begin{array}{cc}
0 & D^*_0\\
D_0 & 0\\
\end{array}\right).$$
Then $(\cA,H,D)$ is a $p$-summable even spectral triple with
respect to the grading on $H$ given by
$$\gamma:=\left(\begin{array}{cc}
1 & 0\\
0 & -1\\
\end{array}\right).$$
The construction of Theorem \ref{S3cross} is correspondingly modified
to yield an {\em odd} spectral triple $(\cA\rtimes_\alpha\Z,
\cH,\cD)$ with $\cH=\ell^2(\Z,H)=H\otimes\ell^2(\Z)$, the
Hilbert space for the regular representation of 
$\cA\rtimes_\alpha\Z$ as before, and with the Dirac operator given by 
$\cD:= D \otimes 1+\gamma\otimes\partial$. Again, if $(\cA,H,D)$ is
$(p,\infty)$-summable then $(\cA\rtimes_\alpha\Z,\cH,\cD)$ is 
$(p+1,\infty)$-summable.

\section{Some motivating examples}\label{examples}

We look here at some simple example that give some motivation 
for introducing spectral triples associated to Mumford curves
and justify why it may be interesting to derive invariants from
these spectral geometries that are more refined than the $C^*$-algebra
$O_A$ itself.

\begin{center}
\begin{figure}
\includegraphics[scale=1.2]{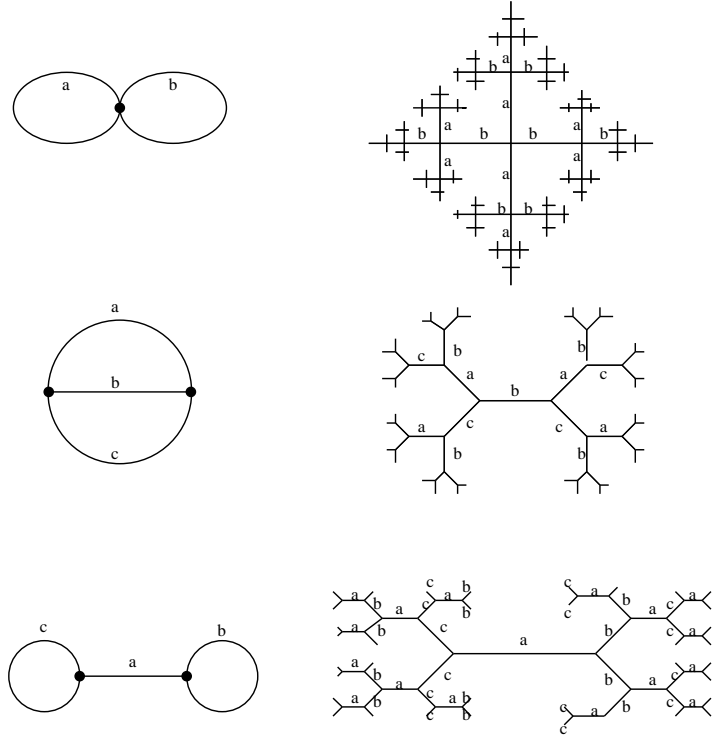}
\caption{The graphs $\cT_\Gamma/\Gamma$
for genus $g=2$, and the corresponding trees $\cT_\Gamma$.
\label{trees}}
\end{figure}
\end{center}

We first look at the case of
genus two Mumford curves discussed in \cite{CM1}, \cite{CM2}. In this
case we are considering a Schottky group $\Gamma$ of rank two in
$\PGL(2,K)$ where $K$ is a finite extension of $\Q_p$. The
combinatorially different forms of the graph $\cT_\Gamma/\Gamma$ of
the special fiber are illustrated in Figure \ref{trees}. 

In terms of corresponding Cuntz--Krieger $C^*$-algebras, we are
considering in this case the algebras $O_{A_i}$, $i=1,2,3$ 
with directed edge matrices
$A_i$ of the form
$$ A_1= \left( \begin{array}{cccc} 1 &1& 0& 1\\
1&1&1& 0  \\
0&1&1& 1 \\
1& 0& 1&1  \end{array} \right) $$
in the first case. In the second case (\cf Figure \ref{trees})
we label by $a=e_1$,
$b=e_2$ and $c=e_3$ the oriented edges in the graph
$\cT_\Gamma/\Gamma$, so that we have a corresponding set of
labels $E=\{ a,b,c,\bar a, \bar b, \bar c \}$ for the edges in the
covering tree $\cT_\Gamma$. A choice of generators for the group
$\Gamma \simeq \Z * \Z$ acting on $\cT_\Gamma$ is obtained by
identifying the generators $g_1$ and $g_2$ of $\Gamma$ with the
chains of edges $a \bar b$ and $a \bar c$. The
directed edge matrix is then of the form
$$ A_2= \left( \begin{array}{cccccc} 0 &1& 0& 0& 0& 1 \\
1&0&1& 0 & 0 & 0 \\
0&1&0& 1& 0& 0 \\
0& 0& 1& 0&1&0 \\
0&0&0&1&0&1 \\
1&0&0&0&1&0 \end{array} \right). $$
The third case in Figure \ref{trees} is analogous. A choice of
generators for the group $\Gamma\simeq \Z * \Z$ acting on
$\Delta_\Gamma$ is given by $ab\bar a$ and $c$.
The directed edge matrix is then
$$ A_3= \left( \begin{array}{cccccc}
0 & 0 & 1 & 0 & 0 & 1 \\
1 & 1 & 0 & 0 & 0 & 0 \\
0 & 0 & 1 & 1 & 0 & 0 \\
0 & 1 & 0 & 0 & 1 & 0 \\
1 & 0 & 0 & 0 & 1 & 0 \\
0 & 0 & 0 & 1 & 0 & 1 \end{array} \right). $$

Recall that, for a Cuntz--Krieger algebra the $K$-theory is
computed in terms of the $n\times n$ matrix $A$ in the form (\cf
\cite{CuKrie}) 
\begin{equation}\label{KOA}
K_0(\overline{O_A})= \Z^n/(1-A^t)\Z^n  \ \ \ \ \  K_1(\overline{O_A})=
{\rm Ker}(1-A^t) \subset \Z^n.
\end{equation}

\smallskip

In general, combinatorially different graphs
$\cT_\Gamma/\Gamma$ are not distinguished by the associated graph 
$C^*$-algebras alone. We can see this in the genus two case as follows. 

\begin{lem}\label{isoOA}
In the case of genus $g=2$, the $C^*$-algebras $\overline{O_{A_i}}$,
$i=1,2,3$, associated to the graphs of Figure \ref{trees} are isomorphic.
\end{lem}

\proof First an explicit calculation shows that the $K$-groups are of
the form 
\begin{equation}\label{KjOAi}
K_j(O_{A_i})\cong \Z^2, 
\end{equation}
for $j=0,1$ and $i=1,2,3$. Moreover, by \cite{Rordam},
for simple Cuntz--Krieger
algebras (\ie algebras $O_A$ where the matrix $A$ is irreducible and
not a permutation matrix) the condition that 
the groups $K_0(O_{A_i})$ are isomorphic implies that the algebras
$\overline{O_{A_i}}$ are isomorphic.
\endproof

\smallskip

Thus, a first question is whether more refined invariants coming
from a spectral triple may be able to distinguish combinatorially
different geometries. There is a more subtle kind of question of
a similar nature. 

As we discussed in \S \ref{TreesSect} above, while the finite graph
$\cT_\Gamma/\Gamma$ only carries the combinatorial information on the
special fiber of the Mumford curve, one can consider the finite graph
$H(\Lambda_\Gamma)/\Gamma$, where $H(\Lambda_\Gamma)$ is the smallest
subtree in the Bruhat--Tits tree $\cT$ of the field $K$ that contains
axes of all elements of $\Gamma$. When one considers the tree 
$H(\Lambda_\Gamma)$ instead of
$\cT_\Gamma$ one is typically adding extra vertices. The way the
tree $H(\Lambda_\Gamma)$ sits inside the Bruhat--Tits tree $\cT$ 
depends on where the Schottky group $\Gamma$ lies in
$\PGL(2,K)$, unlike the information on the graph
$\cT_\Gamma/\Gamma$ which is purely combinatorial (\cf \eg \cite{Ma}).

\smallskip 

\begin{center}
\begin{figure}
\includegraphics[scale=0.4]{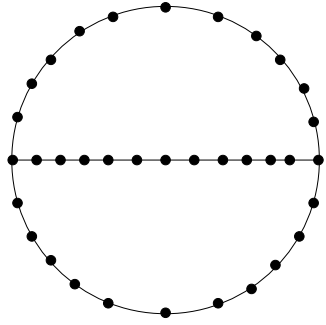}
\caption{The graph $H(\Lambda_\Gamma)/\Gamma$ for $t=\pi^5$ in Kato's
family. \label{graph32}}
\end{figure}
\end{center}

We can consider a specific geometric example, again in the genus two
case, by looking at a 1-parameter family of Mumford curves considered
by Fumiharu Kato in \cite{Kato}. There one considers the free
amalgamated product $N=\Z_m \ast \Z_n$ of two cyclic groups.
For $k$ a (discretely) non-archimedean valued field of characteristic 
coprime to $m$ and $n$, one considers the 
discrete embedding of $N$ in $\PGL(2,k)$  given by
\begin{equation}\label{embedN}
\langle \left( \begin{array}{cc} \zeta^m t & 0 \\ \zeta^m-1 & t
\end{array} \right) \ , \ \left( \begin{array}{cc} \zeta^n &
-(\zeta^n-1)t^{-1}  \\ 0 & 1 \end{array} \right) \rangle . 
\end{equation}
One then considers the free subgroup $\Gamma$ of $N$ generated by
commutators 
\begin{equation}\label{GammaZmn}
 \Gamma := [\Z_m, \Z_n] .
\end{equation}
This group $\Gamma$ is a maximal free subgroup of $N$ of free rank
$g:=(m-1)(n-1)$. It is the Schottky group of
a curve $X_\Gamma$. If one looks at the particularly simple case 
with $m=2$ and $n=3$ one finds a curve of genus two with
$\Z_6$-symmetry. It is not hard to see that the graph of the special
fiber is the second graph in Figure \ref{trees}. Moreover, one can see
that, if $\pi$ denotes a uniformizer and we take $t=\pi^r$ then the
graph $H(\Lambda_\Gamma)/\Gamma$ is again topologically of the same
form but with $2r+1$ vertices inserted on each of the three lines,
namely it looks like the graph in Figure \ref{graph32}.

A direct calculation shows that in this case again the corresponding
algebras $O_{A_r}$ associated to these graphs are stably
isomorphic. In fact, one can show as in the previous example that they
have the same $K$-theory groups. 

\smallskip

This again raises the question of whether more refined invariants, such
as the spectral triple which introduces a smooth subalgebra and
a Dirac operator, can capture more interesting information on the
geometry and the uniformization parameters. 
We hope to return to this general question in future work by
looking more closely at such invariants. For the moment we can give a
heuristic justification of why we think this may be possible. 

\smallskip

In the case of the $\theta$-summable case the information on the
uniformization is stored in the Patterson--Sullivan measure which is
used to define both the representation of the algebra on the Hilbert
space and the Dirac operator. In the case of the finitely summable
construction the dependence on the uniformization will enter again
through a choice of a state $\tau$ on the AF algebra induced by the
the Patterson--Sullivan measure.

\smallskip

The reason why the Patterson--Sullivan measure on the limit set 
is especially good in order to detect geometric properties that depend
on the Schottky uniformization lies in a very general type of rigidity
result (see \cite{Yue}). This type of result implies, in our case,
given two Schottky uniformizations of a fixed genus, if the abstract 
isomorphism of the Schottky groups $\Gamma$ induces a homeomorphism of the
limit sets which is absolutely continuous with respect to the
Patterson--Sullivan measure, then the abstract isomorphism comes from
an automorphism of the ambient group $\PGL(2,K)$. One knows that these
are either inner or they come from automorphisms of the field $K$, so
that one can in fact recover much of the information on the Schottky
uniformization (the Schottky group up to conjugation) from information
on the Patterson--Sullivan measure.

\section{Higher rank cases}\label{higher}

We now proceed to consider some higher dimensional classes of
buildings. The case of rank two are especially interesting because
the classification problem is especially difficult for rank two
and the construction of new invariants can be useful to that
purpose. We first consider the simplest case that can be
reduced to the construction for trees and some simple generalizations 
and then we discuss some more general classes in the hyperbolic case.

\medskip
\subsection{Products of trees}\label{ProdTreesSect}
\hfill\medskip

The first case of 2-dimensional buildings to which the construction
given for trees can be extended is the case
considered in \cite{KiRo}. This deals with affine buildings $\Delta$ whose
2-cells are euclidean squares and whose 1-skeleton is the product of two
trees. On such buildings one considers the action of Burger--Mozes
(BM) groups (\cf \cite{BM}), which are discrete torsion free subgroups
$\Gamma \subset \Aut(\T_1)\times \Aut(\T_2)$ acting freely and
transitively on the set of vertices. If the trees $\T_i$ have even
valences $n_i$, then the quotient $\Delta/\Gamma$ is a polyhedron with one
vertex and $n_1 n_2/4$ square faces, with fundamental group $\Gamma$.
The link of the vertex is a complete bipartite graph with $n_1+n_2$
vertices and edges between each of the first set of $n_1$ vertices and
each of the second set of $n_2$ vertices.

\smallskip

The boundary $\partial \Delta$ is defined by an equivalence relation
on sectors in apartments of $\Delta$, by which two sectors are
equivalent if their intersection contains a sector (this notion
extends the usual shift-tail equivalence of paths in a tree). 
As in the case of trees, the choice of a vertex determines a choice of
a representative in each class of sectors. In the case we are
considering, this gives a non-canonical identification $\partial
\Delta \simeq \partial \T_1 \times \partial \T_2$. 

\smallskip

The action of
$\Gamma$ on $\Delta$ extends to an action on the boundary, hence one can
consider the $C^*$-algebra $C(\partial
\Delta)\rtimes \Gamma$. It is shown in \cite{KiRo} that this is
isomorphic to a rank 2 Cuntz--Krieger algebra associated to two
subshifts of finite type in the alphabet given by $\Gamma$-equivalence
classes of oriented chambers in $\Delta$.
More precisely, the alphabet $\cR$ is given by letters of the form 
$r=(a,b,b',a')$ with $ab=b'a'$, and the horizontal and vertical
transition matrices $A_1$, $A_2$ are commuting $n_1n_2 \times n_1 n_2$
matrices with entries in $\{ 0,1 \}$. They give the admissibility
condition for adjacent squares, in the horizontal and vertical direction,
respectively, \cf \cite{KiRo}. (The notion is well defined, due to the
special class of groups $\Gamma$ considered.)

\smallskip

As in the case of trees, we consider functions on the boundary.
Namely, we consider the Hilbert space $\cL=L^2(\partial\Delta,d\mu)$, where
(after a choice of a base vertex in $\Delta$) the
measure can be identified with a product $d\mu(x) =d\mu_1(x_1) \times
d\mu_2(x_2)$ of Patterson--Sullivan measures on $\partial \T_1$ and
$\partial \T_2$. 

\smallskip

Moreover, the subshifts of finite type in the horizontal and vertical
directions associated to the matrices $A_i$ determine a filtration on
functions on the boundary, analogous to the one considered in the case
of a single tree. In fact, for a fixed base vertex $v_0$ of $\Delta$,
a point $x\in \partial\Delta$ is identified with a sector originating
at $v_0$. Under the identification $\partial\Delta\simeq \partial
\T_1\times \partial \T_2$ determined by $v_0$, we can write, for $f\in
C(\partial\Delta,\Z)\otimes \C$,
$$ f(x)=f(x_1,x_2)=f(a_0a_1a_2\ldots,b_0b_1b_2\ldots) =f(a_0\ldots
a_\ell, b_0\ldots b_k), $$
where $a_0a_1a_2\ldots$ is a path in $\T_1$ and $b_0b_1b_2\ldots$ is a
path in $\T_2$, such that the corresponding infinite 2-dimensional
word is admissible according to the conditions given by the matrices
$A_i$. This means that the space $C(\partial\Delta,\Z)\otimes 
\C$ has a filtration by finite dimensional linear subspaces
$\cV_{\ell,k}$ of functions of finite admissible 2-dimensional words
in the alphabet $\cR$.

\smallskip

We associate a grading operator to this filtration as follows. We 
denote by $\Pi_{\ell,k}$ the orthogonal projection of $\cL$
onto the finite dimensional subspace $\cV_{\ell,k}$. 
We consider the densely defined unbounded self-adjoint grading
operator on $\cL$ of the form 
\begin{equation}\label{inclex}
D = \sum_{m=0}^\infty m \, \sum_{k=0}^m \hat\Pi_{m-k,k}, 
\end{equation}
where $\hat\Pi_{\ell,k}$ is defined by the inclusion-exclusion
$\hat\Pi_{\ell,k}=\Pi_{\ell,k}- \Pi_{\ell-1,k}- \Pi_{\ell,k-1} +
\Pi_{\ell-1,k-1}$. 

\smallskip

\begin{prop}\label{Sp3prodtrees}
Let the data $(\cA,\cH,\cD)$ be given by the algebra
$\cA=C(\partial\Delta)\rtimes \Gamma$ acting on the Hilbert space
$\cH=\cL\oplus \cL$ through a representation \eqref{repU}, and the
operator $\cD=F D$, with $D$ as in \eqref{inclex} and $F$ the
involution exchanging the two copies of $\cL$. 
The data $(\cA,\cH,\cD)$ define a $\theta$-summable spectral
triple. 
\end{prop}

\proof The setting is very similar to the previous cases. We obtain a
representation of the algebra $\cA$ on the Hilbert space $\cH$ by
letting functions in $C(\partial\Delta)$ act on $\cL$ as
multiplication operators and elements in the group $\Gamma$ by
unitary operators
\begin{equation}\label{TopProd}
(T_{\gamma^{-1}}\, f)(x) = \gamma_1'(x_1)^{\delta_1/2} \,
\gamma_2'(x_2)^{\delta_2/2} \, f(\gamma(x)),
\end{equation}
where $x=(x_1,x_2)\in \partial\Delta\simeq \partial
\T_1\times \partial \T_2$ and $\gamma=(\gamma_1,\gamma_2)\in
\Gamma\subset \Aut(\T_1)\times\Aut(\T_2)$. 

\smallskip

In order to show that there exists a dense involutive subalgebra
$\cA_0\subset \cA$, such that the commutators $[\cD,a]$ are bounded
operators on $\cH$ for all $a\in \cA_0$, it is
sufficient to prove that the commutators of $D$ with functions in
$C(\partial\Delta,\Z)\otimes \C$ and with operators $T_\gamma$ for 
$\gamma\in \Gamma$ are bounded. The first case is clear since a
function $h\in C(\partial\Delta,\Z)\otimes \C$ is contained in some 
$\cV_{u,v}$, hence the corresponding multiplication operator maps
$\cV_{\ell,k}$ to itself whenever $\cV_{\ell,k} \supset
\cV_{u,v}$. Since the $\gamma_i'$ are locally constant on $\partial
\T_i$, we also obtain that there exists some integer $n=n(\gamma)$
such that $T_\gamma : \cV_{\ell,k} \to \cV_{\ell+n,k+n}$ whenever
$\cV_{\ell,k} \supset \cV_{u,v}$, for $\gamma_i'\in
\cV_{u,v}$. 
\endproof

\smallskip

It is easy to see that, even if the case considered here is euclidean
and not hyperbolic, still the spectral triple constructed this way
will not be finitely summable.
For example, if the group acting is a product $\Gamma \times \Gamma$
of two copies of a Schottky group of genus $g$ and $\T_1=\T_2$ is the
Cayley graph of $\Gamma$, and $A_1=A_2$ is the directed edge matrix of
the quotient group, namely the matrix with $A_{ij}=1$ unless
$|i-j|=g$, then the dimensions of the eigenspaces of the
operator \eqref{inclex} are $\dim E_m = (m+1) 2g(2g-1)^{m-1}(2g-2)^2$
for $m\geq 2$, $4 g(2g-1)(2g-2)$ for $m=1$ and $2g(2g-1)$ for $m=0$. 

\smallskip

Noncommutative spaces associated to Euclidean buildings, in the form
of crossed product $C^*$-algebras $C(\partial\Delta)\rtimes \Gamma$
were also considered in the case of buildings with triangular
presentations. Of particular arithmetic interest is the case of the
``fake projective planes''. These are algebraic surfaces with ample 
canonical class and the same numbers $p_g=q=0$ and $c_1^2=3c_2=9$ as
the projective plane $\P^2$. The first such example was constructed by
Mumford \cite{Mum2} using p-adic uniformization by a discrete
cocompact subgroup $\Gamma \subset \PGL(3,\Q_2)$. More recent
examples (known as CMSZ fake projective planes) were obtained, again using
actions on the Bruhat--Tits building of $\PGL(3,\Q_2)$ (see \cite{CMSZ},
\cite{KaOc}). The formal model of the CMSZ fake projective planes is
obtained by the action on the Bruhat--Tits building $\Delta$ of
$\PGL(3,\Q_2)$ of two subgroups $\Gamma_1$, $\Gamma_2$ of index 3. 
The quotients $\Gamma_i\backslash \Omega$, where $\Omega$ is
Drinfeld's symmetric space of dimension 2 over $\Q_2$, are
non-isomorphic fake projective planes. All these examples 
admit a description as Shimura varieties (\cf \cite{KaOc}). The
topology of the corresponding noncommutative spaces 
$\cA_i=C(\partial\Delta)\rtimes \Gamma_i$
distinguishes between the fake projective planes. In fact, the results
of \cite{RS3} show for instance that, for the CMSZ cases, the K-theory
of $\cA_i$ is $K_0(\cA_1)=\Z/3$ while $K_0(\cA_2)=\Z/2\oplus\Z/2\oplus\Z/3$.
Recently, Gopal Prasad and Sai-Kee Yeung identified in \cite{PraYeu} the 
complete list of fake projective planes (see also \cite{PraYeu2} for 
a higher dimensional case).  It would be interesting to see
if these can be studied from the operator theoretic point of view and
whether the $K$-theory of the relevant $C^*$-algebras distinguishes
them or whether more refined invariants of spectral geometry
can be used to that purpose.

\smallskip

As we discussed in \S \ref{examples}, one can consider higher 
dimensional cases of  
combinatorially different actions of the same group $\Gamma$ on a
building $\Delta$, such as two
combinatorially non-equivalent presentations
of the same group acting on a product of two trees of degree four
(presentations 42 and 44 from \cite{KiRo}):
$$P_1=\{a_1,b_1,a_2,b_2: a_1b_1a_1^{-1}b_1^{-1}, a_1b_2a_1^{-1}b_2^{-1},
 a_2b_1a_2^{-1}b_2, a_2b_2a_2^{-1}b_1\}$$
$$P_2=\{a_1,b_1,a_2,b_2: a_1b_1a_1^{-1}b_2, a_1b_2a_1^{-1}b_1,
 a_2b_1a_2^{-1}b_2, a_2b_2a_2^{-1}b_1\}.$$
One can ask whether in such examples invariants
coming from spectral triples may be able to distinguish the
combinatorially different actions.

\medskip
\subsection{Polyhedra covered by products of trees}\label{covers}
\hfill\medskip

In the previous section we looked at 2-dimensional buildings that are
products of trees, with the action of groups of BM type, 
$\Gamma\subset \Aut(\T_1)\times \Aut(\T_2)$. In this section we show
that the results of the previous section may be applied more
generally. Namely, 
we show the existence of an infinite family of examples that are not
of BM type, but which can be reduced to BM type by passing to a
subgroup of index four.

\smallskip

A {\em polyhedron} is a two-dimensional
complex obtained from several oriented $p$-gons
by identification of corresponding sides.
Consider a vertex of the polyhedron and 
take a sphere of a small radius at this point.
The intersection of the sphere with the polyhedron is
a graph, which is called the {\em link} at this vertex.

\smallskip

Recall that a graph is {\em bipartite} if its set of vertices
can be partitioned into two disjoint subsets $P$ and $L$ such that no
two vertices in the same subset lie on a common edge.
It is known by the result of \cite{BB} that the universal cover of 
a polyhedron with square faces and complete bipartite graphs 
as links is a 2-dimensional Euclidean buildings which is a product of
two trees $\T_1 \times \T_2$. This means that an efficient method to
construct Euclidean buildings with compact quotients is by
constructing finite polyhedra with appropriate links.

\smallskip

We recall the definition of {\em polygonal presentation} given in
\cite{V}. 

\begin{defn}\label{polygpres}
Suppose given $n$ disjoint connected bipartite graphs
$G_1, G_2, \ldots G_n$.
Let $P_i$ and $L_i$ denote the sets of black and white vertices in
$G_i$, for $i=1,...,n$. Let $P=\cup P_i$ and $L=\cup L_i$, with $P_i
\cap P_j = \emptyset$ and $L_i \cap L_j = \emptyset$, for $i \neq j$. 
Let $\lambda$ be a bijection $\lambda: P\to L$.

A set $\cP$ of $k$-tuples $(x_1,x_2, \ldots, x_k)$, with $x_i \in P$,
is be called a {\em polygonal presentation} over $P$ compatible
with $\lambda$ if the following properties are satisfied.

\begin{enumerate}
\item If $(x_1,x_2,x_3, \ldots ,x_k) \in \cP$, then
   $(x_2,x_3,\ldots,x_k,x_1) \in \cP$.
\item Given $x_1,x_2 \in P$, then  $(x_1,x_2,x_3, \ldots,x_k) \in 
\cP$ for some $x_3,\ldots,x_k$ if and only if $x_2$ and $\lambda(x_1)$
are incident in some $G_i$.
\item Given $x_1,x_2 \in P$, then  $(x_1,x_2,x_3, \ldots ,x_k) 
\in \cP$ for at most one $x_3 \in P$.
\end{enumerate}
If there exists such $\cP$, then the corresponding $\lambda$ is called a
{\em basic bijection}.
\end{defn}

\smallskip

Polygonal presentations for $n=1$, $k=3$, with the incidence graph of
the finite projective plane of order two  
or three as the graph $G_1$, were listed in \cite{CMSZ}.

\smallskip

One can associate a polyhedron $X$ with $n$ vertices to a
polygonal presentation $\cP$ in the following way.
To every cyclic $k$-tuple $(x_1,x_2,x_3,\ldots,x_k)$ 
we assign an oriented $k$-gon, with the word $x_1 x_2 x_3\ldots x_k$
written on its boundary. The polyhedron is then obtained by 
identifying sides with the same labels in these 
polygons, preserving orientation.
We say then that the polyhedron $X$ corresponds to the polygonal
presentation $\cP$.
It was shown in \cite{V} that a polyhedron $X$ that corresponds to
a polygonal presentation $\cP$ has the
graphs $G_1, G_2, \ldots, G_n$ as links.

\smallskip

In particular, suppose that $X$ is a polyhedron corresponding to
a polygonal presentation $\cP$ and let $s_i$ and $t_i$
be, respectively, the number of vertices and of edges of the 
graph $G_i$, for $i=1,...,n$. Then $X$ has $n$ vertices
(the number of vertices of $X$ is equal to the number of graphs),
$k \sum_{i=1}^n s_i$ edges and $\sum_{i=1}^n t_i$ 
faces. All the faces are polygons with $k$ sides.

\smallskip

We use the procedure illustrated above to construct certain
compact polyhedra with square faces whose links are complete bipartite
graphs. 

\begin{defn}\label{bipart} We say that the sets $\mathcal{G}_1,
\mathcal{G}_2,\dots, \mathcal{G}_k$
of connected bipartite graphs are {\em compatible}, if
all graphs in every set have the same number $n$ of white vertices   
and, for every pair $\mathcal{G}_1$, $\mathcal{G}_j$, $j=2,\dots,k$,
there is a corresponding bijection between sets of white vertices that
preserves the degrees of the vertices.
\end{defn}

\begin{prop}\cite{V4} Let $\mathcal{G}_1,
\mathcal{G}_2,\dots, \mathcal{G}_k$ be compatible sets of
connected bipartite graphs, $k\ge1$. Then there exists a family
of finite polyhedra with $2k$-gonal faces, whose links at the vertices
are isomorphic to the graphs from $\mathcal{G}_1, \mathcal{G}_2,\dots,
\mathcal{G}_k$.
\end{prop}

We give now an explicit construction of a particular case of this
theorem, when $k=2$ and each of the families $\mathcal{G}_i$, $i=1,2$ 
contains exactly one complete bipartite graph $G_i$, with $n$ white
and $r$ black vertices.

\smallskip

By \cite{V}, to construct the polyhedron with given links,
it is sufficient to construct a corresponding polygonal presentation.
By the definition of compatible sets of bipartite graphs,
there is a  bijection $\alpha$ from the set
of white vertices of  $G_1$
to the set of white vertices of $G_2$, preserving the degrees
of the white vertices.
We mark the white vertices of $G_i$, $i=1,2$ by
letters of an alphabet $\mathcal{A}_i=\{ x_1^i, \ldots, x_n^i\}$, 
such that the bijection $\alpha$ is induced by the indices of the letters,
\ie $\alpha(x_m^1)=x_m$.
We mark the black vertices of $G_i$, $i=1,2$ by
letters of an alphabet $\mathcal{B}_i=\{y_1^i, y_2^i, \ldots,
y_r^i\}$. Thus, every edge of $G_i, i=1,2$  can be
presented in a form $(x_m^iy_l^i)$, for $m=1,\dots,n$ and
$l=1,\dots,r$. 

\smallskip

Having such a  bijection $\alpha$ of white
vertices we can choose a bijection $\beta$ of the
set of edges of $G_1$ to the set of edges of
$G_2$, which preserves $\alpha$. Let
$\beta_j(x^1_m y^1_l)=x_m^2 y_j^2$. We let the cyclic
word $(x^1_m, y^1_l,x_m^2, y_j^2)$, for $m=1,\dots,n$ and  
$l,j=1,\dots,r$ belong to
the set $\mathcal{P}$. It is shown in \cite{V4} (in more
general form) that $\mathcal{P}$ is a
polygonal presentation.
Denote then by $X$ the polyhedron corresponding to this polygonal
presentation $\mathcal{P}$.

\begin{defn}\label{stablepair}
A polygonal presentation $\mathcal{P}$
satisfies the {\em stable pairs condition} if any word
$(x^1_m, y^1_l,x_m^2, y_j^2) \in \mathcal{P}$ if and only if 
every word in $\mathcal{P}$ containing $x^1_m$ or $x_m^2$ 
has the form $(x^1_m,y^1_s,x^2_m,y^2_t)$ and
every word in $\mathcal{P}$ containing $y^1_l$ or $y_j^2$ 
has the form $(x^1_p,y^1_l,x^2_p,y^2_j)$.
\end{defn}

With the condition of Definition \ref{stablepair}, we obtain the
following result.

\begin{lem}\label{stpairBM}
If $X$ is a polyhedron $X$ corresponding to a polygonal
presentation that satsifies the stable pairs
condition, then the fundamental group of $X$ is of BM type.
\end{lem}

\proof
Consider a polyhedron $X$ with
polygonal presentation  $\mathcal{P}$=$(x^1_m, y^1_l,x_m^2,
y_j^2)$, for $m=1,\dots,n$ and  
$l,j=1,\dots,r$, which satisfies the stable pair condition.
This polyhedron has four vertices and $nr$ faces with words from 
$\mathcal{P}$ on their boundary. To compute its fundamental group
we need letters from one word of $\mathcal{P}$, say $(x^1_m,
y^1_l,x_m^2, y_j^2)$, to be trivial. Because of the stable pair
condition, $y^1_sy^2_t=1$, for $s,t=1,\ldots,r$, 
and $x^1_px^2_p=1$. All relations of the fundamental group of $X$ have
the form $(x^1_p, y^1_s,(x^1_p)^{-1}, (y^1_s)^{-1})$, for
$p=1,\dots,m-1,m+1,\ldots, n$ and $s=1,\dots,t-1,t+1,\ldots, r$. It is
then a BM group by definition 
and it acts on a product of trees of valences $2(n-1)$ and $2(r-1)$.
\endproof

\smallskip

We now consider the following infinite family of examples.
Let $G$ be a complete bipartite graph on $8q$ vertices,
$4q$ black vertices and $4q$ white ones. Let $\mathcal{A}$ and
$\mathcal{B}$ be two alphabets on $4q$  letters,
$\mathcal{A}=\{x_1, x_2, \ldots, x_{4q}\}$ and $\mathcal{B}=\{y_1, y_2, \ldots,
y_{4q}\}$.
We mark every black vertex with an element from
$\mathcal{A}$ and every white vertex with an element from $\mathcal{B}$.

\smallskip

We define a polygonal presentation $\mathcal{P}$
as the following set of cyclic words
\begin{equation}\label{cyclwords}
\begin{array}{lr}
 (x_{1+4i},x_{2+4j},x_{4+4i},x_{3+4j}), & 
 (x_{1+4i},x_{1+4j},x_{4+4i},x_{4+4j}), \\[3mm]
 (x_{1+4i},x_{3+4j},x_{4+4i},x_{2+4j}),  &
 (x_{2+4i},x_{2+4j},x_{3+4i},x_{3+4j}), 
\end{array}
\end{equation}
where $i,j=0,1, \ldots, q-1$. The basic bijection $\lambda$ is given by
$\lambda(x_l)=y_l$.

\smallskip

The polyhedron $X$ that corresponds to $\mathcal{P}$ has 
square faces and one vertex whose link is
naturally isomorphic to a complete bipartite graph. By \cite{BB} the
universal covering $\Delta$ of $X$ is the direct product of two
trees. 

\smallskip

For example, in the case $q=1$ the polygonal presentation $\mathcal{P}$
contains four cyclic words $(x_{1},x_{2},x_{4},x_{3})$,
$(x_{1},x_{1},x_{4},x_{4})$, $(x_{1},x_{3},x_{4},x_{2})$,
$(x_{2},x_{2},x_{3},x_{3})$. The corresponding polyhedron $X$ 
consists of four faces and one vertex, with the link at this
vertex given by the complete bipartite graph with four vertices of
each color. 

\smallskip

In this class of examples, the fundamental group $\Gamma$ of the
polyhedron acts on the 2-dimensional building $\Delta\simeq \T_1\times
\T_2$ so that $X=\Delta/\Gamma$, but $\Gamma$ is not a
subgroup of $\Aut(\T_1)\times \Aut(\T_2)$. 
To reduce this case to the case of BM groups,
we have to find in $\Gamma$ a subgroup of finite index which is 
of BM type.

\smallskip

\begin{lem}\label{BMfinind} 
The group $\Gamma$ with generators $x_k, k=1,...,4q$
and relations
\begin{equation}\label{rels}
\begin{array}{lr}
 x_{1+4i}x_{2+4j}x_{4+4i}x_{3+4j}=1, & 
 x_{1+4i}x_{1+4j}x_{4+4i}x_{4+4j}=1, \\[3mm]
 x_{1+4i}x_{3+4j}x_{4+4i}x_{2+4j}=1, &  
 x_{2+4i}x_{2+4j}x_{3+4i}x_{3+4j}=1,
\end{array}
\end{equation}
where $i,j=0,1, \ldots, q-1$, contains a subgroup of index four which
is of BM type. 
\end{lem}

\proof
Consider the polyhedron $Y$ which corresponds to the polygonal
presentation \eqref{cyclwords}. It contains one vertex and $4q$ faces
and the group $\Gamma$ is the fundamental group of $Y$.
The 4-branching cover of $Y$ is a polyhedron with four vertices
of type $X$, which satisfies the stable pairs condition. Indeed, the
polygonal presentation of $Y$ can be obtained \eqref{cyclwords} by
replacing each word $(x_k,x_l,x_m,x_n)$ by four:  
$(x_k^1,x_l^2,x_m^3,x_n^4)$,
$(x_k^4,x_l^1,x_m^2,x_n^3)$,
$(x_k^3,x_l^4,x_m^1,x_n^2)$,
$(x_k^2,x_l^3,x_m^4,x_n^1)$.
By direct inspection we can see that this satisfies the stable pair condition.
By Lemma \ref{stpairBM}, the fundamental group
of $X$ is of BM type, acting on the product of two trees of valence
$2(4q-1)$. Thus, the group $\Gamma$ contains
a subgroup $\Gamma_0$ of index four of BM type.
\endproof

\smallskip

In this case, one can still apply the construction of Proposition
\ref{Sp3prodtrees} to the index four subgroup of Lemma
\ref{BMfinind}. This means that we now work with the algebra
$C(\Delta)\rtimes \Gamma_0$, which is Morita equivalent to
$C(\Delta\times S)\rtimes \Gamma$, for $S=\Gamma/\Gamma_0$ the coset
space with the left action of $\Gamma$, instead of working with
$C(\Delta)\rtimes \Gamma$.

\medskip
\subsection{Hyperbolic 2-dimensional buildings}\label{2HypSect}
\hfill\medskip

Finally, we give a construction of a spectral triple analogous
to the one defined on trees, in the case of a class of hyperbolic
2-dimensional buildings. We will concentrate on the {\em right angled
Fuchsian buildings}. We recall briefly some properties of these
buildings (\cf \cite{BP}). They are obtained as universal cover of
orbihedra, according to the following construction. One begins with a
regular $r$-gon $P\subset \H$ with angles $\pi/2$. With the edges
labeled clockwise with $\{ i \}$ ($i=1,\ldots,r$) and vertices
correspondingly labeled $\{ i, i+1\}$ and $\{ r, 1 \}$, for assigned
labels $q_i\geq 2$, one obtains a orbihedron by assigning the trival
grou to the face of $P$, the cyclic group $\Gamma_i=\Z/(q_i+1)\Z$ to the $\{ i
\}$ edge and the group $\Gamma_i\times \Gamma_{i+1}$ to the $\{ i, i+1
\}$ vertex. This orbihedron has universal cover $\Delta$, with link at
an $\{ i \}$ labeled vertex given by the complete bipartite graph on
$(q_i+1)+(q_{i+1}+1)$ vertices. The complex $\Delta$ is a hyperbolic
building (in fact a Tits building), where every apartment is
isomorphic to the hyperbolic plane $\H$ with the tessellation
given by the action on $P$ of the cocompact Fuchsian Coxeter group generated by
inversions on the edges. Recall also that a wall in $\Delta$ is
a doubly infinite geodesic path in the 1-skeleton of $\Delta$. By the
form of the links, all edges in such a path have the same label $\{ i
\}$. The equivalence relation on edges given by being in the same wall
has the tree-walls as equivalence classes. A tree-wall with label $\{ i
\}$ divides $\Delta$ into $(q_i+1)$ components. The boundary
$\partial\Delta$ is defined by the set of geodesic rays from a base
point in $\Delta$. It is homeomorphic to the universal Menger curve
$\cM$. Isometries of $\Delta$ extend to homeomorphisms of the
boundary. 

\smallskip

For a discrete finitely generated $\Gamma \subset {\rm Isom}(\Delta)$,
the limit set $\Lambda_\Gamma \subset \partial\Delta$ is
the set of accumulation points of orbits of $\Gamma$. The group
$\Gamma$ is nonelementary if $\Lambda_\Gamma$ consists of more than
two points. We can also consider, as in the case of trees, the geodesic
hull $H(\Lambda_\Gamma)$, obtained by considering all infinite geodesics in
$\Delta$ with endpoints on the limit set $\Lambda_\Gamma$. The group
$\Gamma$ is quasi-convex-cocompact if $H(\Lambda_\Gamma)/\Gamma$ is
compact (\cf \cite{Coo}). 
For instance, we can consider the case where $\Gamma$ is the
fundamental group of the orbihedron and $\partial
\Delta=\Lambda_\Gamma$.  

\smallskip

If $\cG_\Delta$ denotes the dual graph of $\Delta$, with a vertex for
each chamber and an edge connecting two chambers whenever these 
share an edge in $\Delta$. If such edge is of type $\{ i \}$ the
corresponding edge in $\cG_\Delta$ is given length $\log q_i$. This
defines a metric $d(v,w)$ on $\cG_\Delta$. The horospherical distance
is then defined as in the case of trees \eqref{horodist},
\eqref{gammaprime}, by setting $(v_1,v_2,v)=d(v_1,v)-d(v_2,v)$ which
induces a well defined function $(v_1,v_2,x)$, for $x\in
\partial\Delta$, satisfying the cocycle relation
$(v_1,v_2,x)=(v_1,v,x)+(v,v_2,x)$. The function $(v_1,v_2,x)$ is
locally constant on $\partial_{reg}\Delta$, the complement in 
$\partial\Delta$ of the set of endpoints of walls of $\Delta$. 

\smallskip

On $\partial\Delta$ one considers the combinatorial metric
$\delta_v(x,y)$ (\cf \cite{BP}). If $\delta_H$ is the Hausdorff
dimension and $d\mu_v$ the Hausdorff measure, for $\gamma\in {\rm
Isom}(\Delta)$ one has 
\begin{equation}\label{NxvwMeas}
(\gamma^* d\mu_v)(x) = e^{\delta_H \tau (v,\gamma^{-1}v,x)} \, d\mu_v(x),
\end{equation}
where $\tau$ is the unique positive solutions to the equation (\cf
\cite{BP}) 
\begin{equation}\label{eqtau}
\sum_{i=1}^r \frac{q_i^x + q_{i+1}^x}{(1+q_i^x)(1+q_{i+1}^x)} =2.
\end{equation}

\smallskip

Consider a covering $\{ U_{n,i} \}_{n\geq 1, 1\leq
i\leq i_n}$ of the Menger curve $\cM\simeq \partial\Delta$, such that 
$U_{n,i}=\partial\Delta\cap C_{n,i}$, where as a cell complex each
$C_{n,i}$ is a cube and the set of vertices of the $C_{n,i}$ form a
$\Gamma$-invariant subset of $\partial\Delta\smallsetminus
\partial_{reg}\Delta$. For $\Omega$ the closure of
$\partial_{reg}\Delta$, consider the $C^*$-algebra
generated by the characteristic functions of
$\Omega \cap U_{n,i}$. This is a commutative 
$C^*$-algebra $\cB$, which is isomorphic to the algebra  
of continuous functions on a
space $\widehat{\partial\Delta}$, which we call the {\em
disconnection} of $\partial\Delta$, by analogy to the spaces
considered in \cite{Spi} in the case of Fuchsian groups acting on
$\P^1(\R)$. The function 
\begin{equation}\label{gammaHyp}
\gamma^\prime_v(x):=e^{\delta_H \tau (v,\gamma^{-1}v,x)}
\end{equation}
is a locally constant function of $x\in \widehat{\partial\Delta}$.

\smallskip

Since $\partial_{reg}\Delta$ is of full $\mu_v$-measure in
$\partial\Delta$, we can identify the Hilbert spaces
$\cL=L^2(\partial\Delta,d\mu_v)$ and $L^2(\widehat{\partial\Delta},d\mu_v)$.
There is on $\cL$ a filtration $\cV$, where $\cV_n$ is the span of
the characteristic functions $\chi_{ \Omega\cap
U_{n,i} }$, for $i=1,\ldots,i_n$. We can consider the corresponding
grading operator $D=\sum_n n \hat\Pi_n$, with $\Pi_n$ the orthogonal
projection onto $\cV_n$ and $\hat\Pi_n=\Pi_n-\Pi_{n-1}$.

\smallskip

The algebra $\cB$ acts on $\cL$ by
multiplication operator, and for $\Gamma$ the fundamental group of the
orbihedron, we have an action by unitaries
\begin{equation}\label{Gammaacthyp}
(T_{\gamma^{-1}}\, f) (x) = \gamma'_v(x)^{1/2} \, f(\gamma x),
\end{equation}
with $\gamma'_v$ as in \eqref{gammaHyp}. This gives a representation
on $\cL$ of the crossed product $\cA=\cB\rtimes \Gamma$.

\begin{thm}\label{2d3S}
Consider the data $(\cA,\cH,\cD)$, where $\cA=\cB\rtimes \Gamma$, as
above, $\cH=\cL\oplus \cL$ and $\cD=FD$, with $F$ the involution
exchanging the two copies of $\cL$. Let $\cA_0\subset \cA$ be the dense
subalgebra generated algebraically by the characteristic functions
$\chi_{ \Omega\cap U_{n,i} }$ and by the elements of $\Gamma$.
For $U$ an automorphism of $\cA$ preserving $\cA_0$, consider
the representation \eqref{repU} of $\cA$ in $\cB(\cH)$. The data
obtained this way determine a $\theta$-summable spectral triple.
\end{thm}

\proof We check the condition on the commutators with $\cD$. The
dense involutive subalgebra of $\cA_0\subset \cA$ is given by the
algebraic crossed product $\cA_0 =\cB_{alg}\rtimes \Gamma$, where
$\cB_{alg}$ is the subalgebra of $\cB$ generated algebraically by the
$\chi_{ \Omega\cap U_{n,i} }$. If $U$ preserves the subalgebra $\cA_0$, it is
sufficient to show that the commutators $[D,a]$ are bounded for all
$a\in \cA_0$. This is clear for $a\in \cB_{alg}$ since such
elements are in some $\cV_k$ for some $k\geq 0$, hence the
corresponding multiplication operators map $\cV_n$ to $\cV_n$ for all
$n\geq k$. In the case of group elements, the commutator
$[T_\gamma,D]$ is also bounded, since $\gamma'_v(x)^{1/2}$ is locally
constant hence, for some $k=k(\gamma)$, we have $T_\gamma: \cV_n \to
\cV_{n+k(\gamma)}$ for all sufficiently large $n$.

\endproof

\medskip
\subsection{Finite summability in higher rank}\label{fintehigher}
\hfill\medskip

There is a setting similar to the one considered in \S \ref{finsum3}
in the higher rank case. In fact,
in \cite{RS}, Robertson and Steger considered
affine buildings $\Delta$ of type $\tilde A_2$, whose boundary $\Lambda$ is
defined by an equivalence relation on sectors (just as in the case of
trees it is given by an equivalence relations on geodesics). They
showed that, if $\Gamma$ is a group of type rotating automorphisms of
$\Delta$, then the $C^*$-algebra $C(\Lambda)\rtimes \Gamma$ is
isomorphic to a higher rank Cuntz--Krieger algebra $O_{A_1,A_2}$.
This is a particular (rank two) case of more general higher rank
generalizations of Cuntz--Krieger 
algebras, associated to a finite collection of transition matrices
$A_j$, $j=1,\ldots, r$, with entries in $\{ 0, 1 \}$, associated to
shifts in $r$ different directions, with the transition matrices
satisfying compatibility conditions (see conditions (H0)--(H3) of
\cite{RS}). The matrices give admissibility conditions for
$r$-dimensional words in an assigned alphabet.
In the case of the Cuntz--Krieger algebra
$O_A=C(\Lambda_\Gamma)\rtimes \Gamma$ one can choose as 
generators the partial isometries $S_{u,v}=T_{uv^{-1}}P_v$, for
$u,v\in \Gamma$, with $t(u)=t(v)$ (same tail as edges in the Cayley
graph). Similarly, in the higher rank case, one has generators
that are partial isometries $S_{u,v}$, where $u$ and $v$ are words in
the given alphabet, with $t(u)=t(v)$. These satisfy the relations
\begin{equation}\label{highCKrel}
\begin{array}{ll}
S_{u,v}^*=S_{v,u} &\ \  S_{u,v}S_{v,w}=S_{u,w} \\[2mm]
S_{u,v}=\sum S_{uw,vw} &\ \  S_{u,u}S_{v,v}=0, \,\, \forall u\neq v
\end{array}
\end{equation}
The sum here is over $r$-dimensional words $w$ with source
$s(w)=t(u)=t(v)$ and with shape $\sigma(w)=e_j$, for $j=1,\ldots,r$, where $e_j$
is the $j$-th standard basis vector in $\Z^r$ (see \cite{Ro},
\cite{RS} for more details). 
Robertson and Steger also proved (\cite{RS} \S 6) that the higher rank
Cuntz--Krieger algebras $O_{A_1,\ldots, A_r}$ are stably isomorphic
to a crossed product
\begin{equation} \label{AF-T-high} \overline{O_{A_1,\ldots, A_r}} \cong
\overline{\mathcal F}_{A_1,\ldots, A_r} \rtimes_T \Z^r, \end{equation}
where $\cF_{A_1,\ldots, A_r}$ is the AF algebra generated by the
$S_{u,v}$ with $\sigma(u)=\sigma(v)$. 
Again the stabilization $\overline{\cF_A}$ is a non-unital AF algebra.
One expects that a similar technique, based on the standard
spectral triple of the $n$-torus $T^n$ and a spectral triple for
the non-unital AF algebra $\overline{\mathcal F}_{A_1,\ldots, A_r}$
to yield finitely summable triples for the algebra of
\eqref{AF-T-high}.


\medskip

\begin{thebibliography}{99} 

\bibitem{AntChris} C.~Antonescu, E.~Christensen, {\em Spectral triples
for AF $C^*$-algebras and metrics on the Cantor set}, preprint
arXiv math.OA/0309044.

\bibitem{sB83}
S.~Baaj, P.~Julg, {\em
Th\'eorie bivariante de Kasparov et op\'erateurs non born\'es dans les 
$C^*$-modules hilbertiens}, C. R. Acad. Sci. Paris S\'er. I Math. 
296 (1983), no. 21, 875--878.

\bibitem{BB} W.~Ballmann, M.~Brin, {\em Orbihedra of nonpositive
curvature}, Publications Math\'ematiques IHES, 82 (1995), 169--209.

\bibitem{BP} M.~Bourdon, H.~Pajot, {\em Rigidity of quasi--isometries
for some hyperbolic buildings}, Comment. Math. Helv. 75 (2000)
701--736. 

\bibitem{BM} M.~Burger, S.Mozes, {\em Lattices in products of trees},
Publ. Math. IHES, 92 (2001), 151--194.

\bibitem{CMSZ}  D.I.~Cartwright, A.M.~Mantero, T.~Steger, A.~Zappa, 
{\em Groups acting simply transitively on the vertices of a 
building of type $\tilde A\sb 2$. I/II}.  
Geom. Dedicata  47  (1993), 143--166 and 167--223. 

\bibitem{Connes2} A.~Connes, {\em Compact metric spaces, Fredholm
modules, and hyperfiniteness}, Ergodic Theory and Dynamical Systems 9
(1989) 207--220.

\bibitem{Co94}  A.~Connes, {\em Noncommutative geometry}, Academic
Press, 1994.

\bibitem{Connes} A.~Connes, {\em Geometry from the spectral point of
view}. Lett. Math. Phys. 34 (1995), no. 3, 203--238.

\bibitem{CoMo} A.~Connes, H.~Moscovici, {\em The local index formula
in noncommutative geometry}.  Geom. Funct. Anal.  5  (1995),  no. 2,
174--243. 

\bibitem{CM} C.~Consani, M.~Marcolli, {\em Noncommutative geometry,
dynamics and $\infty$-adic Arakelov geometry},  
Selecta Math. (N.S.)  10  (2004),  no. 2, 167--251.  

\bibitem{CM1} C.~Consani, M.~Marcolli, {\em Spectral triples from
Mumford curves}, International Math. Research Notices, 36 (2003)
1945--1972.

\bibitem{CM2} C.~Consani, M.~Marcolli, {\em New perspectives in
Arakelov geometry},  Number theory,  81--102, CRM Proc. Lecture Notes,
36, Amer. Math. Soc., Providence, RI, 2004. 

\bibitem{Coo} M.~Coornaert, {\em Mesures de Patterson--Sullivan sur le
bord d'un espace hyperbolique au sens de Gromov}, Pacific J. Math. 159
(1993) 241--270.

\bibitem{CoKa} G.~Cornelissen, F.~Kato, {\em Equivariant deformation
of Mumford curves and of ordinary curves in positive characteristic},
Duke Math. J. 116 (2003) 431--470.

\bibitem{CuKrie} J.~Cuntz, W.~Krieger, {\it A class of $C^*$--algebras
and topological Markov chains, I,II}, Invent. Math. 56 (1980) 251--268
and Invent. Math. 63 (1981) 25--40.

\bibitem{HeHu} S.~Hersonsky, J.~Hubbard, {\em Groups of automorphisms
of trees and their limit sets}, Ergod. Th. Dynam. Sys. 17 (1997)
869--884. 

\bibitem{Moyal} V.~Gayral, J.M.~Gracia-Bondia, B.~Iochum,
T.~Sch\"ucker, J.C.~Varilly, {\em Moyal Planes are Spectral Triples}, 
Comm. Math. Phys.  246  (2004),  no. 3, 569--623. 

\bibitem{Kato} F.~Kato, {\em Non-Archimedean orbifolds covered by
Mumford curves}.  J. Algebraic Geom. Vol.14  (2005),  no. 1, 1--34. 

\bibitem{KaOc} F.~Kato, H.~Ochiai, {\em Arithmetic structure of CMSZ
fake projective planes}, preprint.

\bibitem{KiRo} J.S.~Kimberley, G.~Robertson, {\em Groups acting on
products of trees}, New York J. Math. 8 (2002) 111--131.

\bibitem{LiMa} D.~Lind, B~Marcus, {\em An introduction to symbolic
dynamics and coding}. Cambridge University Press, Cambridge, 1995.  

\bibitem{Lu} A.~Lubotsky, {\em Lattices in rank one Lie groups over
local fields}, Geometric and Functional Analysis, 1 (1991) 405--431.

\bibitem{Man} Yu.I.~Manin, {\em  Three-dimensional hyperbolic geometry
as $\infty$-adic Arakelov geometry},  Invent. Math.  104  (1991),
no. 2, 223--243. 

\bibitem{Ma} Yu.I.~Manin, {\em $p$-adic automorphic
functions}. Journ. of Soviet Math., 5 (1976) 279-333.

\bibitem{Matsu} K.~Matsumoto, {\em On automorphisms of $C^*$-algebras
associated with subshifts}, J. Operator Theory, Vol.44 (2000) 91--112.

\bibitem{Mum} D.~Mumford, {\em An analytic construction of
degenerating curves over complete local rings}, Compositio Math.
24 (1972) 129--174.

\bibitem{Mum2} D.~Mumford, {\em An algebraic surface with $K$ ample,
$K^2=9$, $p_g=q=0$}, Amer. J. Math., Vol. 101, (1979) N.1, 233--244.

\bibitem{PraYeu} G.~Prasad, S.K.~Yeung, {\em Fake projective planes},
math.AG/0512115.

\bibitem{PraYeu2} G.~Prasad, S.K.~Yeung, {\em Cocompact arithmetic 
subgroups of $PU(n-1,1)$ with Euler-Poincare characteristic $n$ and 
a fake $P^4$}, math.AG/0602144. 

\bibitem{Ro} G.~Robertson, {\em Boundary actions for affine buildings
and higher rank Cuntz--Krieger algebras}, in ``$C\sp *$-algebras'' 
(M\"unster, 1999),  182--202, Springer, 2000. 

\bibitem{RS} G.~Robertson, T.~Steger, {\em Affine buildings, tiling
systems and higher rank Cuntz-Krieger algebras}.  J. Reine
Angew. Math.  513  (1999), 115--144.  

\bibitem{RS3} G.~Robertson, T.~Steger, {\em $K$-theory computations
for boundary algebras of $\tilde A_2$ groups}, preprint.

\bibitem{Rordam} M.~Rordam, {\em Classification of Cuntz-Krieger 
algebras}. $K$-Theory  9  (1995), no.1, 31--58.

\bibitem{Spi} J.S.~Spielberg, {\em Cuntz--Krieger algebras associated
with Fuchsian groups}, Ergod. Th. Dynam. Sys. Vol. 13 (1993) 581--595. 

\bibitem{Sull} D.~Sullivan, {\em On the ergodic theory at infinity of
an arbitrary discrete group of hyperbolic motions}.
Riemann surfaces and related topics: Proceedings of the 1978 Stony
Brook Conference (State Univ. New York, Stony
Brook, N.Y., 1978), pp. 465--496, Ann. of Math. Stud., 97, Princeton
Univ. Press, 1981.

\bibitem{V} A.~Vdovina, 
{\em Combinatorial structure of some hyperbolic buildings},
Math. Z. 241 (2002), no. 3, 471--478.

\bibitem{V2}
A.~Vdovina, {\em Polyhedra with specified links},
S\'emin. Th\'eorie Spectr. G\'eom., Institut Fourier,
Universit\'e Grenoble I, 21-37, (2003).

\bibitem{V4} 
A.~Vdovina, {\em Groups, periodic planes and buildings},
J. Group Theory  8  (2005),  no. 6, 755--765.

\bibitem{Wer} A.~Werner, {\em  Arakelov intersection indices of linear
cycles and the geometry of buildings and symmetric spaces},  Duke
Math. J.  111  (2002),  no. 2, 319--355.

\bibitem{Yue} C.~Yue, {\em Mostow rigidity of rank $1$ discrete groups with 
ergodic Bowen-Margulis measure.}
Invent. Math. 125 (1996), no. 1, 75--102.

\end{thebibliography}
\end{document}